\documentclass[preprint,11pt]{elsarticle}

\newtheorem{theorem}{Theorem}[section]
\newtheorem{lemma}{Lemma}[section]
\newtheorem{corollary}{Corollary}[section]

\newtheorem{remark}{Remark}[section]
\newcommand{\ignore}[1]{}{}

\usepackage{amssymb}
\usepackage{amsmath}

\usepackage{color}
\usepackage{soul}
\usepackage[dvipsnames]{xcolor}

\usepackage[colorlinks=true, linkcolor=blue, citecolor=blue]{hyperref}

\def\1{{{\mbox{${\rm{1\negthinspace\negthinspace I}}$}}}}

\renewcommand{\theequation}{\arabic{equation}}

\newcommand\beq{\begin{equation}}
\newcommand\eeq{\end{equation}}
\renewcommand{\theequation}{\thesection.\arabic{equation}}

\usepackage{ifthen}
\usepackage{xkeyval}
\usepackage{todonotes}
\begin{document}

\begin{frontmatter}

\title{Comparison on the criticality parameters for two supercritical branching processes in random environments}
\author[cor1]{Xiequan Fan$^{*,}$}
\author[cor2]{Haijuan Hu}
\author[cor3]{Hao Wu}
\author[cor4]{Yinna Ye}
 \cortext[cor1]{\noindent Corresponding author. \\
\mbox{\ \ \ \ }\textit{E-mail}:  fanxiequan@hotmail.com (X. Fan)\\
\mbox{\ \ \ \ \ \ \ \ \ \ \  \ \ \ \ } . }
\address[cor1,cor2]{School of Mathematics and Statistics, Northeastern University at Qinhuangdao, Qinhuangdao, China}
\address[cor3]{Center for Applied Mathematics,
Tianjin University, Tianjin 300072, China}
\address[cor4]{Department of Applied Mathematics, School of Science, Xi'an Jiaotong-Liverpool University, 215123 Suzhou, China  }

\begin{abstract}
Let $\{Z_{1,n} , n\geq 0\}$ and $\{Z_{2,n}, n\geq 0\}$ be two  supercritical branching processes in different random environments,
with criticality parameters $\mu_1$ and $\mu_2$ respectively. It is known that with certain conditions,  $\frac{1}{n} \ln  Z_{1,n}  \rightarrow \mu_1$
and $\frac{1}{m} \ln  Z_{2,m}  \rightarrow \mu_2$ in probability as $m, n \rightarrow \infty.$ In this paper, we are interested in the comparison on the two criticality parameters, which  can be regarded as two-sample $U$-statistic.   To this end, we prove a  non-uniform Berry-Esseen's bound  and    Cram\'{e}r's  moderate deviations for $\frac{1}{n} \ln  Z_{1,n}   -  \frac{1}{m} \ln  Z_{2,m}$ as $m, n \rightarrow \infty.$
An application is also given for constructing confidence intervals
 of $\mu_1-\mu_2$.
\end{abstract}

\begin{keyword} Branching processes; Random environment;  Berry-Esseen's bound; Cram\'{e}r's  moderate deviations
\vspace{0.3cm}
\MSC primary 60F10; 60K37; secondary 60J80
\end{keyword}

\end{frontmatter}




\section{Introduction}
\setcounter{equation}{0}

The branching process in a random environment (BPRE) was first  introduced by Smith and Wilkinson \cite{SW69} to model the growth of a population in an independent and identically distributed (iid) random environment.
 Various limit theorems for BPRE  have been obtained: basic results for extinction probabilities and limit theorems for  BPRE can be found in Athreya and Karlin \cite{AK71a,AK71b}.
For  subcritical BPRE,  researches focus on the study of the survival probability and conditional limit theorems: see,
for instance, Vatutin \cite{V10}, Afanasyev \emph{et al.}\ \cite{ABV14}, Vatutin and Zheng \cite{VZ12}  and  Bansaye  and  Vatutin \cite{BV17}.
While, for supercritical  BPRE,   a number of researches have studied  limit theorems, see, for instance,
 Gao \emph{et al.}\ \cite{GQW2014}, Hong and Zhang \cite{HZ19}, Gao \cite{G2021} and
Li \emph{et al.}\ \cite{LHP22}. For moderate and large deviations, we refer to
  B\"{o}inghoff  and Kersting \cite{BK10}, Bansaye and  B\"{o}inghoff \cite{BB11},
  Huang and Liu \cite{HL12},  Nakashima \cite{N13}, B\"{o}inghoff \cite{B14},
   Grama \emph{et al.}\ \cite{GLE17} and   \cite{FHL20}.
See also Wang and Liu \cite{WL17} and Huang \emph{et al.}\ \cite{HWW22} for BPRE with immigrations.

Let's introduce, first of all, the classical model of BPRE. Suppose $\xi=\left(\xi_{n}\right)_{n\geq0}$ is a sequence of iid random variables. Usually, $\xi$ is called environment process and is related to the environment condition of the system, which may be captured by a random process. Consider a particle system placed in a random environment  $\xi$. In the system, the population size of particles at each generation is denoted by a discrete-time random process $(Z_n)_{n\geq0}$, which evolves as follows.  Given an environment process $\xi$, at generation 0, assume that there is only one particle $\emptyset$ and $Z_0=1$; at  generation  1, the particle $\emptyset$ splits and is replaced by a group of $Z_1$ new particles; and those new particles constitute the population in the first generation; at  generation  2, every individual, say the individual $i$, in the generation $1$ splits independently with the others and is replaced by a group of  $X_{1,i}$ new particles and  the new particles produced from
the individuals in the  first generation constitute the population in the second generation, and so forth. More generally, in the $n$-th ($n\geq1$) generation, there are $Z_n$ particles in the system and in the $(n+1)$-th generation, each particle, say $i$, splits independently and is replaced by a group of $X_{n,i}$ new particles. In mathematical words, a discrete-time random process $(Z_n)_{n\geq0}$ is called BPRE in the environment  $\xi$, if it satisfies the following recursive relation:
$$
Z_{0}=1, \ Z_{n+1}=\sum_{i=1}^{Z_{n}} X_{n, i}, \ n \geq 0,
$$
where $X_{n, i}$ represents the number of  particles produced by the $i$-th individual in the $n$-th generation.  The distribution of $ X_{n,i}$, which is also called the offspring distribution, depending on the environment $ \xi_{n} $ is denoted by $p\left(\xi_{n}\right)=\left\{p_{k}\left(\xi_{n}\right)=\mathbb{P}(X_{n,i}=k |\xi_{n}): k \in \mathbb{N}\right\}$. Suppose that given $ \xi_{n} $, $ (X_{n, i})_{i\geqslant1}$ is a sequence of iid random variables; moreover, $ (X_{n, i})_{i\geqslant1}$ is independent of $(Z_{1},\ldots ,Z_{n}) $. Let $\left(\Gamma , \mathbb{P}_{\xi}\right)$ be the probability space under which the process is defined when the environment $\xi$ is given. The state space of the environment process $ \xi $ is denoted by $ \Theta $ and the total probability space can be regarded as the product space $\left(\Theta^{\mathbb{N}}\times\Gamma, \mathbb{P}\right),$ where $\mathbb{P}(d x, d \xi)=\mathbb{P}_{\xi}(d x) \tau(d \xi).$ That is, for any measurable positive function $g $ defined on $\Theta^{\mathbb{N}}\times\Gamma$, we have

$$\int g(x, \xi) \mathbb{P}(d x, d \xi)=\iint g(x, \xi) \mathbb{P}_{\xi}(d x) \tau(d \xi),$$
where $ \tau $ is the distribution law of the environment process $ \xi $. And $\mathbb{P}_{\xi} $ can be regarded as the conditional probability when the environment $ \xi $ is given. Usually,  the conditional probability $\mathbb{P}_{\xi}$ is called quenched law, while the total probability  $\mathbb{P}$  is  called  annealed law.

In this paper, we are interested in two branching processes placed separately in two different environments $\xi_1$ and $ \xi_2$. More precisely, let  $(\xi_1, \xi_2)^T :=( (\xi_{1,n}, \xi_{2,n})^T)_{n\geq 0} $ be a sequence of iid random vectors, where $T$  is the transport operator and $ ( \xi_{1,n}, \xi_{2,n})^T \in \mathbb{R}^2$ stands for the random environment at generation $n$.
Thus $(\xi_{1,n}, \xi_{2,n})^T_{n\geq 0}$ are independent random vectors, but for given $n$, $\xi_{1,n}$ and $\xi_{2,n}$ may not be independent.
For any $n \in \mathbb{N}$, each realization of $\xi_{1,n}$ corresponds to a probability law $ p(\xi_{1,n} ) = \{ p_i(\xi_{1,n}): i \in  \mathbb{N}\},$   that is $p_i(\xi_{1,n}) \geq 0  $ and $ \sum_{i=0}^\infty p_{i}(\xi_{1,n})=1$. Similarly, each realization of $\xi_{2,n}$ corresponds to a probability law $ p(\xi_{2,n})$.
Let $\{Z_{1,n}\}_{n\geq 0}$ and $\{Z_{2,n}\}_{ n\geq 0}$ be two  branching processes in the random environment
$\xi_1$ and $\xi_2$, respectively.
Then $\{Z_{1,n}\}_{n\geq 0}$ and $\{Z_{2,n}\}_{ n\geq 0}$  can be described  as follows: for $n \geq 0,$
$$Z_{1,0}=1,\quad   Z_{1,n+1} = \sum_{i=1}^{Z_{1,n}} X_{1, n,i},\quad  Z_{2,0}=1,\quad  Z_{2,n+1} = \sum_{i=1}^{Z_{2,n}} X_{2,n,i} ,  $$
where $X_{1,n,i}$ and $X_{2,n,i}$ are the number of  particles produced by  the $i$-th individual in generation $n$ in the environments $\xi_{1}$ and  $\xi_{2}$,
respectively. Moreover,  we assume that
given $(\xi_{1,n}, \xi_{2,n})^T$, the random variables $\{X_{1, n, i},  X_{2,n,i}, i\geq 1 \}$   are independent, and they are also independent of $\{Z_{1,k}, Z_{2,k},  0 \leq k \leq n\}$.
 Denote $\mathbb{P}_{\xi_1 ,\xi_2}$  the conditional probability when the environment
 $(\xi_1,\, \xi_2)^T$ is given, and  $\tau $  the joint distribution law  of the environment  $(\xi_1,\, \xi_2)^T$. Then
 $$\mathbb{P}(dx_1, dx_2,  dy_1, dy_2)=\mathbb{P}_{\xi_1, \xi_2} (dx_1, dx_2)\tau (dy_1, dy_2)$$
is the joint  annealed law of the two branching processes placed separately in two different environments.

 In particular, if   $\xi_1$ and $\xi_2$ are  independent, then we have $ \tau (dy_1, dy_2)= \tau_1 (dy_1 ) \tau_2 (  dy_2)$, where $\tau_1$ and $\tau_2$ are the marginal distributions of $\xi_1$ and $\xi_2$ respectively. In the sequel,  $\mathbb{E}_{\xi_1, \xi_2}$ and $\mathbb{E}$ denote the expectations with respect to $\mathbb{P}_{\xi_1, \xi_2}$  and $\mathbb{P}$, respectively. For any $ n\geq 1$,
set
$$m_{1,n}^{(p)}  = \sum_{k=0}^{\infty} k^p \, p_k(\xi_{1,n} ),  \quad m_{2,n}^{(p)}   = \sum_{k=0}^{\infty} k^p \, p_k(\xi_{2,n}  ),    $$
$$ \Pi_{1,n} = \prod_{i=0}^{n-1} m_{1,i}, \quad \Pi_{2,n}  = \prod_{i=0}^{n-1} m_{2,i},  $$
with the convention that $ \Pi_{1,0} = \Pi_{2,0} = 1$. Clearly,  $\{m_{1,n}^{(p)}\}_{n\geq 0}$ and $\{m_{2,n}^{(p)}\}_{n\geq 0}$ are
two sequences of iid random variables.
 For simplicity of notations,  write
  $$m_{1,n}= m_{1,n}^{(1)} \ \ \  \ \ \   \textrm{and} \ \ \  \ \ m_{2,n}=m_{2,n}^{(1)} ,$$
and denote
 \begin{align*}
 M_{1,n} =  \ln m_{1,n} ,\quad  M_{2,n} =  \ln m_{2,n},  \quad
 \mu_1  = \mathbb{E}M_{1,0}, \quad \mu_2  = \mathbb{E}M_{2,0}, \\
 \sigma_1 ^2 =  \textrm{Var}(M_{1,0} ) , \quad  \sigma_2 ^2  =  \textrm{Var}( M_{2,0}),\quad \displaystyle   \ \ \ \ \ \ \ \  \ \ \ \
 \end{align*}
where $\mu_1$  and $\mu_2$ are known as  the criticality parameters for
BPREs $\{Z_{1,n}\}_{ n\geq 0}$ and $\{Z_{2,n}\}_{n\geq 0}$, respectively.
 To avoid  the environments $\xi_1$ and $\xi_2$ are degenerated,  assume that $  0 <\sigma_1,  \sigma_2  < \infty.$
 Denote
  \begin{align*}
   \rho  =  \frac{\textrm{Cov}(M_{1,0}, M_{2,0})}{\sigma_1\sigma_2 }
 \end{align*}
the correlation coefficient between $M_{1,0}$ and $M_{2,0}$. In particular, if $\xi_1$ and $\xi_2$ are independent, we have $\rho=0.$ Write  $ \ln ^+  x  = \max\{ \ln x , 0 \}. $
 Throughout the paper,  assume the following conditions  used in Grama \emph{et al.}\ \cite{GLE17}:
\begin{equation}\label{dgf103hm}
\mathbb{E}\bigg[ \frac{Z_{1,1}}{m_{1,0}} \ln^+ Z_{1,1}   +  \frac{Z_{2,1}}{m_{2,0}} \ln^+ Z_{2,1} \bigg]  < \infty
\end{equation}
and
\begin{equation} \label{dgf1sm}
p_0(\xi_{1,0})=p_0(\xi_{2,0} ) =0, \ \ \  \ \ \ \ a.s.
\end{equation}
The assumption (\ref{dgf1sm}) implies  that each individual has at least one offspring.  The assumptions  (\ref{dgf103hm}) and (\ref{dgf1sm})
together imply that the processes $\{Z_{1,n}\}_{ n\geq 0}$ and $\{Z_{2,n}\}_{ n\geq 0}$ are both supercritical, and satisfy that $\mu_1, \mu_2 >0$  and $\mathbb{P}(Z_{1, n} \rightarrow \infty )= \mathbb{P}(Z_{2, n} \rightarrow \infty )=1$. See Athreya and Karlin \cite{AK71b} and Tanny \cite{T88}.

For the case of a single supercritical BPRE, say $\{Z_{1,n}\}_{ n\geq 0}$, the normal approximation  has been well studied.
With the additional conditions $ \mathbb{E}  \big(\frac{Z_{1,1}}{m_{1,0}} \big )^{p} < \infty$ for a constant $ p>1 $ and   $\mathbb{E}M_{1,0}^{2+\rho}    < \infty$
  for a constant $\rho \in (0, 1)$,  Grama et al.\ \cite{GLE17}   have established the following Berry-Esseen bound  for $\ln Z_{1,n}$:
\begin{equation} \label{bercrser}
\sup_{x\in \mathbb{R}} \Big|\mathbb{P}\Big( \frac{ \ln  Z_{1,n}  - n \mu_1 \ }{  \sigma_1 \sqrt{n}} \leq  x  \Big)  -  \Phi(x)   \Big|
\leq  \frac{  C }{   n^{\rho/2}   }.
\end{equation}
Assume
$ \mathbb{E}  \frac{Z_{1,1}^{p}}{m_{1,0}}   < \infty$ for a constant $ p>1$ and    $\mathbb{E}e^{\lambda_0 M_{1,0}}   < \infty$
for a constant $\lambda_0>0$,
  Grama et al.\ \cite{GLE17} have also  established the following Cram\'{e}r moderate deviation expansion:
for $0 \leq x = o(\sqrt{n} )$ as $n\rightarrow \infty$,
\begin{equation}\label{cramer}
\Bigg| \ln \frac{\mathbb{P}\Big( \frac{ \ln  Z_{1,n}  - n \mu_1 \ }{  \sigma_1 \sqrt{n}} \geq x  \Big)}{1-\Phi(x)} \Bigg|  \leq C    \frac{  1+x^3    }{  \sqrt{n}  },
\end{equation}
where  $C$ is a positive  constant. Such type of inequality  implies the following result about the equivalence to the normal tail, that is
\begin{equation}\label{cramedr}
\frac{\mathbb{P}\big( \frac{ \ln  Z_{1,n}  - n \mu_1 \ }{  \sigma_1 \sqrt{n}}   \geq x  \big)}{1-\Phi(x)} =1+o(1)
\end{equation}
uniformly for  $x \in [0, o(n^{1/6}))$  as $n\rightarrow \infty$.
The results (\ref{bercrser})-(\ref{cramedr}) are interesting both in theory and in applications. For example,  when the parameter $\sigma_1$ is known, they can be applied to construct  confidence intervals  for estimating the criticality parameter $\mu_1$ in terms of the observation $Z_{1,n}$ and the generation $n$.

Despite the fact that the limit theorems for one supercritical BPRE are well studied, there is no result
for comparison on the criticality parameters  for two supercritical BPREs.
The objective of the paper is to fit up this gap.
Consider the following common hypothesis testing:
$$H_0: \mu_1-\mu_2 = 0\ \     \ \ \ \textrm{versus} \ \ \  \ \ \ H_1: \mu_1 -\mu_2\neq 0.$$
When $\mu_1$ and $\mu_2$ are  means of two independent populations, such type of hypothesis testing has been considered  by Chang \emph{et al}.\ \cite{CSZ16}, and Cram\'{e}r type moderate deviations has been established therein.
 In this paper, we are interested in the case $ \mu_1$  and $\mu_2$ are two criticality parameters of BPREs.
 Notice that $\ln  Z_{i,n} =\sum_{k=1}^{n}M_{i,k}+ \ln (Z_{i,n}/\Pi_{i,n}),  i=1, 2.  $
 With some mild conditions, we have $\frac 1n\ln (Z_{i,n}/\Pi_{i,n}) \rightarrow 0$ in probability as $  n \rightarrow \infty$.
As $\sum_{k=1}^{n}M_{i,k}$, $i=1, 2$, are both sums of iid random variables, by the law of large numbers,  we have  $\frac{1}{n} \ln  Z_{1,n}  \rightarrow \mu_1$ and $\frac{1}{m} \ln  Z_{2,m}  \rightarrow \mu_2$ in probability as $m\wedge n \rightarrow \infty$.
Thus, to test the hypothesis,
 it is crucial to estimate  the asymptotic distribution of the random variable  $\frac{1}{n} \ln  Z_{1,n} - \frac{1}{m} \ln  Z_{2,m}$,
 which is the main  purpose of this paper.
Notice that $\frac{1}{n} \ln  Z_{1,n} - \frac{1}{m} \ln  Z_{2,m}$ has an asymptotic distribution
 as $\frac{1}{n} \sum_{k=1}^{n}M_{1,k} - \frac{1}{m} \sum_{k=1}^{m}M_{2,k}.$
  When $\xi_1$ and $\xi_2$ are independent,   $\sum_{k=1}^{n}M_{1,k}$ and $\sum_{k=1}^{n}M_{2,k}$  are both sums of iid random variables and
 then the hypothesis testing can be regarded as  a two-sample $U$-statistic problem.

Define
\begin{equation*}
R_{m,n} = \frac{ \frac{1}{n} \ln  Z_{1,n}  -  \frac{1}{m} \ln  Z_{2,m} - (\mu_1 - \mu_2) }{  \sqrt{\frac1n  \sigma_1^2  +  \frac1m  \sigma_2^2 -2 \rho \sigma_1 \sigma_2 \frac{m \wedge n}{m \, n} \ } },\ \ \quad n, m \in \mathbb{N}.
\end{equation*}
Throughout  the paper,  we  also assume   either $$ \rho \in [-1, 1) \ \ \ \ \ \ \  \textrm{or}  \ \ \ \ \ \ \  \rho=1 \ \ \textrm{but} \   \sigma_1   \neq  \sigma_2.$$
The last condition ensures  that $\frac1n  \sigma_1^2  +  \frac1m  \sigma_2^2 -2 \rho \sigma_1 \sigma_2 \frac{m \wedge n}{m \, n}$ is in order of $\frac{1 }{m\wedge n}$ as $m \wedge n \rightarrow \infty$. Indeed, if $m \leq n,$
 it is easy to see that
 $$\frac1n  \sigma_1^2  +  \frac1m  \sigma_2^2 -2 \rho \sigma_1 \sigma_2 \frac{m \wedge n}{m \, n}=(\frac1m -\frac1n)\sigma_2^2 +\frac{ \sigma_1^2  -2 \rho \sigma_1 \sigma_2+\sigma_2^2}{n} \asymp  \frac{1}{m  }  .$$
  Let us introduce our main results briefly. Firstly,  Theorem \ref{th00} presents the central limit theorem (CLT) for $R_{m,n}:$
for all $x \in \mathbb{R},$ it holds
\begin{equation}
   \lim_{m \wedge n \rightarrow \infty }\mathbb{P}\big(  R_{m,n} \leq  x  \big)  = \Phi(x),
\end{equation}
where $\Phi(x)$  is the standard normal distribution function.
Secondly, under the moment condition $\mathbb{E} [ M_{1,0}^{2+\delta}+  M_{2,0}^{2+\delta} \,]  < \infty $ with
$\delta \in (0, 1] $ and the condition $\mathbb{E}\Big[ \frac{Z_{1,1}  ^{p}}{m_{1,0}^p}   +  \frac{Z_{2,1}  ^{p}}{m_{2,0}^p}   \Big]<\infty $ with $p>1$, Theorem \ref{th02} gives a non-uniform  Berry-Esseen bound  for $R_{m,n}$:   for any constant $\delta' \in (0, \delta)$ and all $x \in \mathbb{R},$
\begin{equation}\label{nonfsds}
  \bigg|\mathbb{P}\big(  R_{m,n} \leq  x  \big)  -  \Phi(x)   \bigg|
\leq  \frac{  C }{   (m\wedge n)^{\delta/2}   }\frac{  1 }{  1+|x|^{1+\delta'}  }.
\end{equation}
By Lemma  \ref{lemma4.3} and (\ref{4.8}) in the paper, it seems that $R_{m,n}$ only has a finite moment of order $1+\delta'$, $\ \delta' \in (0, \delta)$, under the stated conditions, which explains why the  non-uniform  Berry-Esseen bound  is of order $\displaystyle   |x|^{-1-\delta'}   $ instead of order $\displaystyle    |x|^{-2-\delta}  $ as  $x \rightarrow \infty$.
In particular, when $m \rightarrow \infty,$ we have $\frac{1}{m} \ln  Z_{2,m} \rightarrow \mu_2$ in probability, which leads to $R_{m,n} \rightarrow  \frac{ \ln  Z_{1,n}  - n \mu_1 \ }{  \sigma_1 \sqrt{n}}$ in  probability. Thus   inequality (\ref{nonfsds}) implies that
\begin{equation}\nonumber
  \bigg|\mathbb{P}\Big( \frac{ \ln  Z_{1,n}  - n \mu_1 \ }{  \sigma_1 \sqrt{n}} \leq  x  \Big)  -  \Phi(x)   \bigg|
\leq  \frac{  C }{   (m\wedge n)^{\delta/2}   }\frac{  1 }{  1+|x|^{1+\delta'}  },
\end{equation}
which  improves the  Berry-Esseen bound (\ref{bercrser}) by adding a factor $\frac{  1 }{  1+|x|^{1+\delta'}  }$. Moreover, using
this non-uniform  Berry-Esseen bound, we can obtain an optimal convergence rate in the Wasserstein-$1$ distance.
 Thirdly, we  establish Cram\'{e}r's moderate deviations.
Assuming Cram\'{e}r's condition  $\mathbb{E} \big[ e^{\lambda_0 M_{1,0} } +   e^{\lambda_0 M_{2,0} }  \big]  < \infty $
for a constant $\lambda_0>0$ and $\mathbb{E} \Big[ \frac{ Z_{1,1}  ^{p}}{m_{1,0}} +   \frac{ Z_{2,1} ^{p}}{m_{2,0}}  \Big] <\infty$ for a constant $p>1$,  Theorem \ref{th2.2}   shows   that   for all $0 \leq x   \leq C^{-1} \sqrt{m \wedge n }  $,
\begin{equation}   \label{tsh02s}
 \Bigg| \ln \frac{\mathbb{P}\big( R_{m,n} \geq x  \big)}{1-\Phi(x)} \Bigg|
  \leq C      \frac{ 1+x^3   }{  \sqrt{m \wedge n} \ }.
\end{equation}
 From  \eqref{tsh02s}, we  obtain
the following equivalence to the normal tail: it holds
\begin{equation} \label{fgdsx}
\frac{\mathbb{P}\big(    R_{m,n} \geq x  \big)}{1-\Phi(x)} =1+o(1)
\end{equation}
uniformly  for   $0\leq x= o((m \wedge n)^{1/6}) $  as $m\wedge n\rightarrow \infty$.
When  $m\rightarrow \infty$, it is easy to see
that (\ref{tsh02s}) and (\ref{fgdsx}) remain valid when $R_{m,n}$ is  replaced by $\frac{ \ln  Z_{1,n}  - n \mu_1 \ }{  \sigma_1 \sqrt{n}} $.
Thus our results recover Cram\'{e}r's moderate deviations (\ref{cramer}) and (\ref{cramedr}) established by
Grama \emph{et al.}\ \cite{GLE17}. Finally,
as an application  of our results, we discuss  the construction of confidence intervals for $\mu_1-\mu_2$.

The paper is organized as follows. Our main results are stated and discussed in Section \ref{sec2}.
 In Section 3, an  application  of our results to construction of confidence intervals for $\mu_1-\mu_2$  is demonstrated.
The proofs of the main results are given in Section \ref{sec4d}.

Throughout the paper, $c$ and $C,$ probably supplied with some indices,
denote respectively a small positive  constant and a large positive constant.
Their values may vary from line to line. For two sequences of positive  numbers $(a_n)_{n\geq 1}$ and $(b_n)_{n\geq 1}$,
we write $ a_n  \asymp b_n $ if there exists a positive constant $C$ such that for all $n,$ it holds $C^{-1}b_n \leq a_n \leq C b_n$.

\section{Main results} \label{sec2}
\setcounter{equation}{0}

In the sequel, we will need the following   conditions.
 \begin{description}
\item[A1]  There exists a constant $\delta \in (0, 1]$
 such that
\[
\mathbb{E} [ M_{1,0}^{2+\delta}+  M_{2,0}^{2+\delta} \,]  < \infty.
\]
\item[A2] There exists a constant $p>1$   such that
$$
\mathbb{E}\bigg[ \frac{Z_{1,1}  ^{p}}{m_{1,0}^p}   +  \frac{Z_{2,1}  ^{p}}{m_{2,0}^p}   \bigg]<\infty .
$$
\end{description}
Denote $\Phi(x)$  the standard normal distribution function.
Let
 \begin{equation*}
 \displaystyle  V_{m,n,\rho} =  \sqrt{\frac1n  \sigma_1^2  +  \frac1m  \sigma_2^2 -2 \rho \sigma_1 \sigma_2 \frac{m \wedge n}{m \, n}}\,,
 \end{equation*}
and recall the definition
\begin{equation} \label{def-Zn0n}
R_{m,n} = \frac{ \frac{1}{n} \ln  Z_{1,n} -  \frac{1}{m} \ln  Z_{2,m} - (\mu_1  - \mu_2) }{ V_{m,n,\rho}}, \quad n, m \in \mathbb{N}.
\end{equation}
We have the following CLT for $R_{m,n}$.
\begin{theorem}\label{th00}
For all $x \in \mathbb{R},$ it holds
\begin{equation}\nonumber
   \lim_{m \wedge n \rightarrow \infty }\mathbb{P}\big(  R_{m,n} \leq  x  \big)  = \Phi(x).
\end{equation}
\end{theorem}

The following theorem gives a non-uniform  Berry-Esseen bound  for $R_{m,n}$.
\begin{theorem}\label{th02}
Assume that the conditions   \textbf{A1} and \textbf{A2} are satisfied.  Let $\delta'$ be a constant such that $\delta' \in (0, \delta).$ Then for all $x \in \mathbb{R},$
\begin{equation}\label{Berr01}
  \Big|\mathbb{P}\big(  R_{m,n} \leq  x  \big)  -  \Phi(x)   \Big|
\leq  \frac{  C }{   (m\wedge n)^{\delta/2}   }\frac{  1 }{  1+|x|^{1+\delta'}  }.
\end{equation}
\end{theorem}

\begin{remark}
By \eqref{4.8} and Lemma  \ref{lemma4.3}, under the conditions \textbf{A1} and \textbf{A2},
$R_{m,n}$ has a finite moment of order $1+\delta'$, $\ \delta' \in (0, \delta)$, which explains why the
non-uniform  Berry-Esseen bound \eqref{Berr01}
is of order $\displaystyle   |x|^{-1-\delta'}   $ instead of order $\displaystyle    |x|^{-2-\delta}  $ as  $x \rightarrow \infty$.
\end{remark}

The following corollary is a direct consequence of Theorem \ref{th02}, which gives a convergence rate of  $R_{m,n}$ to the standard normal random variable in the Wasserstein-$1$ distance.  We first recall the definition of the Wasserstein-$1$ distance.  The Wasserstein-$1$ distance between two distributions $\mu$ and  $\nu$  is defined as follows:
\begin{eqnarray*}
	W_1(\mu, \nu)= \sup  \bigg \{ |\mathbb{E}[f(X)] -\mathbb{E}[f(Y)]    |:   \ (X, Y) \in \mathcal{L}(\mu, \nu),\ f\  \textrm{is  $1$-Lipschitz}     \bigg\} ,
\end{eqnarray*}
where  $\mathcal{L}(\mu, \nu)$ is the collection of all pairs of random variables whose marginal distributions are $\mu$ and $\nu$ respectively.
In particular, if $\mu_M$ is  the distribution of a random variable $M$  and $\nu$ is the standard normal distribution, then
we have   $$W_1(\mu_M , \nu )= d_{w}\left(M\right):=\int_{-\infty}^{+\infty}\Big|\mathbb{P}\left(M \leq x\right)-\Phi(x)\Big| dx , $$
see R\"ollin \cite{AR18}.

\begin{corollary}
\emph{Assume that the conditions   \textbf{A1} and \textbf{A2} are satisfied.  Then}
$$ d_{w}\left(R_{m,n}\right)\leq  \frac{  C }{   (m\wedge n)^{\delta/2}   } . $$
\end{corollary}

By Theorem \ref{th02}, we also have the following Berry-Esseen bounds for $R_{m,n}$.
\begin{corollary} \label{co02}
\emph{Assume that the conditions   \textbf{A1} and \textbf{A2} are satisfied.  Then}
\begin{equation}
\sup_{x\in \mathbb{R}} \Big|\mathbb{P}\big(  R_{m,n} \leq  x  \big)  -  \Phi(x)   \Big|
\leq  \frac{  C }{   (m\wedge n)^{\delta/2}   }  .
\end{equation}
\end{corollary}

Notice that $\frac{1}{m} \ln  Z_{2,m}   \rightarrow  \mu_2  $ in probability as $m\rightarrow \infty$, and thus $$R_{\infty,n} :=\lim_{m \rightarrow \infty }R_{m,n} =   \displaystyle \frac{  \ln  Z_{1,n}  - n \mu_1  }{ \sigma_1 \sqrt{  n  } } $$ in probability. Then  when $m\rightarrow \infty,$ Corollary \ref{co02} recovers  the Berry-Esseen bound established by Grama \emph{et al.}\ \cite{GLE17}, that is,
\begin{equation}
\sup_{x\in \mathbb{R}} \bigg|\mathbb{P}\Big(  \frac{ \ln  Z_{1,n}  - n \mu_1 \ }{  \sigma_1 \sqrt{n}}  \leq  x  \Big)  -  \Phi(x)   \bigg|
\leq  \frac{  C }{   n^{\delta/2}   }  . \nonumber
\end{equation}
It is known that the convergence rate of the last Berry-Esseen bound coincides the best possible one for iid random variables with finite
moments of order $2+\delta$.

Next, we are going to establish  Cram\'{e}r's  moderate deviations for $R_{m,n}$.
To this end, we need the following conditions.
\begin{description}
\item[A3]  The random variables $M_{1,0} $ and $M_{2,0} $ have exponential moments, i.e.\ there exists a constant $\lambda_0>0$
 such that
\[
   \mathbb{E} \big[ e^{\lambda_0 M_{1,0} } +   e^{\lambda_0 M_{2,0} }  \big]    < \infty.
\]

\item[A4]  There exists a constant $p>1$ such that
$$
   \mathbb{E} \bigg[ \frac{ Z_{1,1}  ^{p}}{m_{1,0}} +   \frac{ Z_{2,1} ^{p}}{m_{2,0}}  \bigg] <\infty.
$$
\end{description}
By the definition of $M_{1,0}$ and  $M_{2,0}$,   condition \textbf{A3} is equivalent to $ \mathbb{E} \big[   m_{1,0} ^{\lambda_0  } +    m_{2,0}  ^{\lambda_0  }  \big]    < \infty.$

We have the following  Cram\'{e}r's  moderate deviations for $R_{m,n}$.
\begin{theorem}\label{th2.2}
Assume that the conditions   \textbf{A3} and \textbf{A4} are satisfied.
  Then
 for all $0 \leq   x \leq  c\, \sqrt{m \wedge n}  ,$
\begin{equation}\label{ineq02}
\Bigg|\ln  \frac{\mathbb{P}\big(  R_{m,n}  \geq x  \big)}{1-\Phi(x)} \Bigg| \leq C  \frac{ 1+ x^3   }{  \sqrt{m \wedge n}\  }  .
\end{equation}
\end{theorem}

Thanks to the symmetry  between $m$ and $n$, Theorems \ref{th00}, \ref{th02} and   \ref{th2.2}  remain valid when $R_{m,n}  $ is replaced by  $ -R_{m,n}  $.

By an argument similar to the proof of Corollary 2.2 in \cite{FGLS19}, it is easy to see that
Theorem  \ref{th2.2}  implies  the following moderate deviation principle (MDP) result for $ R_{m,n}$. 
\begin{corollary}\label{corollary02}
\emph{Assume that the conditions   \textbf{A3} and \textbf{A4} are satisfied.
Let $a_n$ be any sequence of real numbers satisfying $a_n \rightarrow \infty$ and $a_n/\sqrt{m\wedge n}\rightarrow 0$
as $ m\wedge n \rightarrow \infty$.  Then, for each Borel set $B$, the following inequalities hold}
\begin{align}
- \inf_{x \in B^o}\frac{x^2}{2} &\leq  \liminf_{n\rightarrow \infty}\frac{1}{a_n^2} \mathbb{P}\bigg(  \frac{R_{m,n}   }{ a_n  }     \in B \bigg) \nonumber \\
 &\leq \limsup_{n\rightarrow \infty}\frac{1}{a_n^2}\ln  \mathbb{P}\bigg(\frac{ R_{m,n} }{ a_n  }  \in B \bigg) \leq  - \inf_{x \in \overline{B}}\frac{x^2}{2}   ,   \label{MDP}
\end{align}
\emph{where $B^o$ and $\overline{B}$ denote the interior and the closure of $B$, respectively.}
\end{corollary}

%
%
%


From  Theorem \ref{th2.2}, using the inequality $|e^y-1| \leq e^C |y|$ valid for $|y| \leq C,$  we obtain the following result about the uniform equivalence to the normal tail.
\begin{corollary}\label{co02s}
\emph{Assume that the conditions   \textbf{A3} and \textbf{A4} are satisfied.  Then it holds}
\begin{equation} \label{rezero}
\frac{\mathbb{P}\big(  R_{m,n}  \geq x  \big)}{1-\Phi(x)} =1+ o(1)
\end{equation}
\emph{uniformly for $x \in [0, \, o((m \wedge n)^{1/6}))$ as $ m\wedge n \rightarrow \infty$.   The result remains valid when $\frac{\mathbb{P} (   R_{m,n} \geq x  )}{1-\Phi(x)}$ is replaced by  $\frac{\mathbb{P} (  R_{m,n}  \leq- x  )}{ \Phi(-x)}$.}
\end{corollary}


Notice that when $\{Z_{2,n}, n\geq 0\}$ is an independent copy of $\{Z_{1,n}, n\geq 0\}$,   we have $\mu_1 = \mu_2, \sigma_1 = \sigma_2$ and $\rho=0.$ Then, for $m=n,$ it holds $ R_{n,n}= \displaystyle \frac{ \ln  Z_{1,n}  -    \ln  Z_{2,n}    }{ \sqrt{2 n } \sigma_1 }.$
Consequently, by Theorem \ref{th02}, under the conditions   \textbf{A1} and \textbf{A2}, it holds for all $x \in \mathbb{R},$
\begin{equation}\nonumber
  \bigg|\mathbb{P}\bigg( \frac{ \ln  Z_{1,n}  -    \ln  Z_{2,n}    }{ \sqrt{2 n } \sigma_1 } \leq  x  \bigg)  -  \Phi(x)   \bigg|
\leq  \frac{  C }{   n^{\delta/2}   }\frac{  1 }{  1+|x|^{1+\delta'}  }.
\end{equation}
By Corollary \ref{co02s}, under the conditions \textbf{A3} and \textbf{A4},  we have
\begin{equation} \nonumber
\frac{\mathbb{P}\Big(  \frac{ \ln  Z_{1,n}  -    \ln  Z_{2,n}    }{ \sqrt{2 n } \sigma_1 }  \geq x  \Big)}{1-\Phi(x)} =1+ o(1)
\end{equation}
uniformly for $x \in [0, \, o(n^{1/6}))$ as $  n \rightarrow \infty$.

\section{Applications to construction of confidence intervals }\label{sec3}
\setcounter{equation}{0}
In this section, we are interested in constructing confidence intervals for $\mu_1-\mu_2. $
This problem is proposed by Behrens in 1929, and it is also known as Behrens-Fisher's problem because Fisher also has  discussed this problem.
When parameters $ \sigma_1,  \sigma_2,$ and $ \rho$ are known, we can apply Theorems  \ref{th02}  and   \ref{th2.2} to construct confidence intervals for $\mu_1-\mu_2. $

 \begin{theorem}\label{th3.1}
 Let $ \kappa_{m,n} \in (0,1) $. Consider the following two groups of conditions:\\
  \begin{description}
\item[B1]
	 The conditions of Theorem  \ref{th02} hold and
	\begin{equation*}
	\left|\ln \kappa_{m,n}\right|=o\big(\ln (m \wedge n)  \big), \ \ \ \textrm{as}\ m \wedge n\rightarrow \infty .
	\end{equation*}
\item[B2]
	 The conditions of Theorem  \ref{th2.2}  hold and
	 \begin{equation*}
	\left|\ln \kappa_{m,n}\right|=o\big((m \wedge n) ^{1/3} \big), \ \ \ \textrm{as}\ m \wedge n\rightarrow\infty.
	\end{equation*}
  \end{description}
If either $\textbf{B1}$ or $\textbf{B2}$ holds, then  for $n$ large enough, $ \left[A_{m,n},\, B_{m,n}\right] $ is the confidence interval of $ \mu_1-\mu_2 $ with confidence level $ 1-\kappa_{m,n} $, where
	$$A_{m,n}= \frac{1}{n} \ln  Z_{1,n}    -  \frac{1}{m} \ln  Z_{2,m} - V_{m,n,\rho} \Phi^{-1}\left(1-\frac{\kappa_{m,n}}{2}\right)$$
  and
  $$ \quad B_{m,n}=\frac{1}{n} \ln  Z_{1,n}    -  \frac{1}{m} \ln  Z_{2,m} +  V_{m,n,\rho} \Phi^{-1}\left(1-\frac{\kappa_{m,n}}{2}\right).$$	
\end{theorem}
\emph{Proof.}
Assume  that the condition $\textbf{B1}$ is satisfied. By Theorem  \ref{th02}, as $m\wedge n\rightarrow \infty$, we have
\begin{align}\label{k7}
\frac{\mathbb{P}\left(R_{m,n} > x\right)}{1-\Phi(x)}=1+o(1)\quad\text{and}\quad \frac{\mathbb{P}\left(R_{m,n} < -x\right)}{\Phi(-x)}=1+o(1)
\end{align}
uniformly for $ 0\leq x=o\left(\sqrt{\ln (m\wedge n) }\right).$ For $ p\searrow 0 $, the quantile function of the standard normal distribution has the following asymptotic expansion $$ \Phi^{-1}(p)=-\sqrt{\ln \frac{1}{p^{2}}-\ln \ln \frac{1}{p^{2}}-\ln (2 \pi)}+o(1) .$$ In particular,
when $\kappa_{m,n}$ satisfies the condition $\textbf{B1}$, the upper $\left(1-\frac{\kappa_{m,n}}{2} \right)$-th quantile of standard normal distribution satisfies $$  \Phi^{-1}\left(1-\frac{\kappa_{m,n}}{2}\right)=-\Phi^{-1}\left(\frac{\kappa_{m,n}}{2}\right)=O\left(\sqrt{\left|\ln \kappa_{m,n}\right|}\, \right), $$  which is of order $ o\big(\sqrt{\ln(m\wedge n) }\big) .$ Then, applying the last equality to \eqref{k7},  we have as $ m\wedge n\rightarrow \infty $,
\begin{align}\nonumber
\mathbb{P}\left(R_{m, n} > \Phi^{-1}\left(1-\frac{\kappa_{m,n}}{2}\right)\right) \sim \frac{\kappa_{m,n}}{2}
\end{align}
and
\begin{align}\nonumber
\mathbb{P}\left(R_{m, n} < -\Phi^{-1}\left(1-\frac{\kappa_{m,n}}{2}\right)\right) \sim \frac{\kappa_{m,n}}{2}.
\end{align}
Therefore, as $ m \wedge n\rightarrow \infty $,
\begin{align}\nonumber
\mathbb{P}\bigg(-\Phi^{-1}\left(1-\frac{\kappa_{m,n}}{2}\right)\leq R_{m, n} \leq\Phi^{-1}\left(1-\frac{\kappa_{m,n}}{2}\right)\bigg)\sim 1-\kappa_{m,n},
\end{align}
which implies $\mu_1-\mu_2\in[A_{m,n},\,B_{m,n}]$ with probability $1-\kappa_{m,n}$ for $m \wedge n$ large enough.

Now, assume that the condition $\textbf{B2}$ holds. By  Theorem  \ref{th2.2},  as $m\wedge n\rightarrow \infty$,  we have
\begin{align}\label{m13}
	\frac{\mathbb{P}\left(R_{m, n}> x\right)}{1-\Phi(x)}=1+o(1)\quad\text{and}\quad \frac{\mathbb{P}\left(R_{m, n}< -x\right)}{\Phi(-x)}=1+o(1)
\end{align}
uniformly for $ 0\leq x=o( (m\wedge n)^{1/6}).$
When $\kappa_{m,n}$ satisfies the  condition  $\textbf{B2}$, the upper $\left(1-\frac{\kappa_{m,n}}{2} \right)$-th  quantile of the standard normal distribution satisfies $$\Phi^{-1}\left(1-\frac{\kappa_{m,n}}{2} \right)=-\Phi^{-1}\left(\frac{\kappa_{m,n}}{2} \right)=O\left(\sqrt{\left|\ln \kappa_{m,n}\right|}\right),$$ which is of order $o\left((m\wedge n)^{1/6}\right)$. Then by \eqref{m13}, we get, as $m\wedge n\rightarrow \infty$,
\begin{equation*}
	\mathbb{P}\bigg(-\Phi^{-1}\left(1-\frac{\kappa_{m,n}}{2}\right)\leq R_{m, n}\leq\Phi^{-1}\left(1-\frac{\kappa_{m,n}}{2}\right)\bigg)\sim 1-\kappa_{m,n}.
\end{equation*}
From the equality above, the claim of Theorem \ref{th3.1} still holds.
\hfill\qed

When $\{Z_{2,n}, n\geq 0\}$ is an independent copy of $\{Z_{1,n}, n\geq 0\}$,  we can apply Theorems  \ref{th02}  and   \ref{th2.2} to construct confidence intervals for $\sigma_1. $
 \begin{theorem}\label{pr02}
 Let $ \kappa_{n,n} \in (0,1) $.
If either $\textbf{B1}$ or $\textbf{B2}$ holds, then  for $n$ large enough, $  [A_{n},\, B_n] $ is the confidence interval of $ \sigma_1^2$ with confidence level $ 1-\kappa_{n,n} $, where
	$$A_{n}=\frac{ (\ln  Z_{1,n}  -    \ln  Z_{2,n} )^2   }{ 2 n  \chi_{1-\frac12\kappa_{n,n} }^2(1)  } \ \ \ and \ \ \ B_{n}=\frac{ (\ln  Z_{1,n}  -    \ln  Z_{2,n} )^2   }{ 2 n  \chi_{ \frac12\kappa_{n,n} }^2(1)  } $$
with $\chi_{q }^2(1)$ the $q$-quantiles  for chi-squared distribution with one degree of freedom.
\end{theorem}
\emph{Proof.}
Assume that the condition $\textbf{B1}$ holds. By Theorem  \ref{th02}, as $ n\rightarrow \infty$, we have
\begin{align}\label{kds7}
\frac{\mathbb{P}\left(\frac{ (\ln  Z_{1,n}  -    \ln  Z_{2,n} )^2   }{  2 n   \sigma_1^2  } > x\right)}{\mathbb{P}(  \chi^2(1) \geq x)}=1+o(1)
\end{align}
uniformly for $ 0\leq x=o (\sqrt{\ln n } ).$   Then, applying the last equality to \eqref{kds7},  we have, as $ n\rightarrow \infty $,
\begin{align}\nonumber
\mathbb{P}\bigg( \chi_{ \frac12 \kappa_{n,n} }^2(1) \leq  \frac{ (\ln  Z_{1,n}  -    \ln  Z_{2,n} )^2   }{ 2 n  \sigma_1^2  } \leq  \chi_{1-\frac12\kappa_{n,n} }^2(1)\bigg)\sim 1-\kappa_{n,n},
\end{align}
which implies $\sigma_1^2 \in[A_{n},B_n]$ with probability $1-\kappa_{n,n}$ for $  n$ large enough.

Assume that the condition $\textbf{B2}$ holds. We can use the similar argument and complete  the proof of Theorem \ref{pr02}.
\hfill\qed

\section{Proofs of Theorems}\label{sec4d}
\setcounter{equation}{0}

For $l= 1, 2,$ denote the normalized population size  $$W_{l,n}= \frac{Z_{l,n} }{\Pi_{l,n}},\ \ n \geq 0. $$
Then  $(W_{1, n})_{n\geq0}$ and $(W_{2, n})_{n\geq0}$ are both non-negative martingales
 under the annealed law $\mathbb{P}$, with respect to
the natural filtration $$\mathcal{F}_0=\sigma\{\xi_{1}, \xi_{2}\} ,\ \
\ \mathcal{F}_n=\sigma\{\xi_{1}, \xi_{2}, M_{1, k,i},  M_{2, k,i}, 0\leq k \leq n-1, i\geq 1\},\ \  n\geq 1.$$
  Then,  by Doob's martingale convergence theorem, the limit
 $$W_{l, \infty}=\lim_{n\rightarrow \infty} W_{l,n}$$
 exists $\mathbb{P}$-a.s.\ and, by Fatou's lemma, it satisfies $\mathbb{E} W_{l, \infty} \leq 1$.
The conditions \eqref{dgf103hm} and  \eqref{dgf1sm} together imply that
$\mathbb{P}(W_{l,n} >0)=\mathbb{P}(Z_{l,n} \rightarrow \infty)=\lim_{n \rightarrow \infty} \mathbb{P}(Z_{l,n} >0) =1,$
and that the martingale $\{W_{l,n}\}_{n\geq1}$ converges  to $W_{l, \infty}$ in $\mathbb{L}^1(\mathbb{P})$  (see Athreya and Karlin \cite{AK71b} and also Tanny \cite{T88}).

For simplicity of notations,  without loss of generality, we assume that $m \leq n.$
In the sequel,  denote
 $$\eta_{m,n,i}= \frac{M_{1,i-1}-\mu_1}{ n \, V_{m,n,\rho} \  } , \  \ \ i=1,..., n,   \ \ \ \    \textrm{and} \  \ \ \
   \eta_{m,n,n+j}= -\frac{M_{2,j-1}-\mu_2}{m \, V_{m,n,\rho} \  } , \ \ \  j=1,..., m. $$
Then $R_{m,n}$ can be rewritten in the following form
\begin{equation}\label{4.8}
R_{m,n}=\sum_{i=1}^{n+m}  \eta_{ m,n, i} +\frac{\ln W_{1, n}}{ n\, V_{m,n,\rho}  \ } - \frac{\ln W_{2,m}}{ m\, V_{m,n,\rho} \ } .
\end{equation}
Set
$$Y_i = \eta_{ m,n, i} + \eta_{ m,n, n+i} , \ \ \  i=1, ..., m, \ \  \textrm{and} \ \ \ \ Y_i=\eta_{ m,n, i}, \ \ i=m+1, ..., n.  $$
Then $(Y_i)_{1 \leq i \leq n}$ is a finite sequence of  centered and independent random variables, and satisfies $$\sum_{i=1}^n Y_i = \sum_{i=1}^{n+m}  \eta_{ m,n, i} \ \ \ \textrm{ and } \ \ \ \sum_{i=1}^n \mathbb{E}Y_i^2 =1 . $$
Moreover, we have
$$\textrm{Var}(Y_i) \asymp \frac{1}{m}   , \ \ \  i=1, ..., m, \ \  \textrm{and} \ \ \ \ \textrm{Var}(Y_i)\asymp \frac{m}{n^2} , \ \ i=m+1, ..., n,  $$
as $m   \rightarrow \infty.$ From (\ref{4.8}), we get
\begin{equation}\label{4.9}
R_{m,n}=\sum_{i=1}^n Y_i +\frac{\ln W_{1, n}}{ n\, V_{m,n,\rho}  \ } - \frac{\ln W_{2,m}}{ m\, V_{m,n,\rho} \ } .
\end{equation}

\subsection{Proof of Theorem \ref{th00}}
Without loss of generality, we assume that $m \leq n.$
By the CLT for independent random variables, we have $\sum_{i=1}^n Y_i$ converges in distribution to the standard normal random variable
   as $m  \rightarrow \infty$.  Recall that for $l=1,2,$  $W_{l,n}$ converges  to $W_{l, \infty}$ in $\mathbb{L}^1(\mathbb{P})$ as $n \rightarrow \infty$. By the fact that $V_{m,n,\rho}      \asymp   \frac{1  }{\sqrt{m }},$  the random variables $\frac{\ln W_{1, n}}{ n\, V_{m,n,\rho}  \ } $ and $\frac{\ln W_{2,m}}{ m\, V_{m,n,\rho} \ }$   both convergence in probability to $0$ as $m \rightarrow \infty$.
   Hence,    by (\ref{4.9}),  we have  $R_{m,n}$ converges in distribution to the normal  random variable as $m \rightarrow \infty$. This completes the proof of Theorem \ref{th00}.

\subsection{Preliminary Lemmas for Theorem \ref{th02}}

In the proof of Theorem \ref{th02}, we need the following non-uniform Berry-Esseen bound of  Bikelis \cite{B66}. See also Chen and Shao \cite{CS01}
for more general results.
\begin{lemma}\label{lemma4.1}
		Let $(Y_{i})_{1\leq i \leq n}$ be  independent random variables satisfying $\mathbb{E}Y_{i}=0$ and  $\mathbb{E}\left|Y_{i}\right|^{2+\delta}<\infty$  for some positive constant $\delta \in(0,1]$ and all $1\leq i \leq n$. Assume that $\sum\limits_{i=1}^{n}\mathbb{E}Y_{i}^{2}=1$. Then for all $x \in \mathbb{R},$
		$$\left|\mathbb{P}\bigg(\sum_{i=1}^{n}Y_{i}\leq x\bigg)-\Phi(x)\right|\leq \frac{C}{ 1+|x| ^{2+\delta}}\sum\limits_{i=1}^{n} \mathbb{E}\left|Y_{i}\right|^{2+\delta}.$$
\end{lemma}

Consider the  Laplace  transforms of $W_{1, \infty}$ and $W_{2, \infty}$ as follows:  for all $t\geq 0$,
$$\phi_{i,\xi}(t)=\mathbb{E}_{\xi}e^{-tW_{i, \infty}}\quad\text{and}\quad\phi_i (t)=\mathbb{E}\phi_{i, \xi}(t)=\mathbb{E} e^{-tW_{i, \infty}}, \ \ \ i=1, 2.$$
Clearly, as $W_{i, \infty}  \geq 0$ $\mathbb{P}$-a.s., we have $\phi_{i}(t) \in (0, 1]$, $i=1, 2$. Moreover,
we have the following bounds for $\phi_{i}(t), i=1, 2$, as $t \rightarrow \infty.$
 \begin{lemma}\label{lemma4}
		Assume that the conditions \textbf{A1} and \textbf{A2} are satisfied.  Then for $i=1,2,$ it holds
for all $t>0,$
$$\phi_{i}(t) \leq \frac{C}{1+ (\ln^+ t
			)^{1+\delta}}.$$
\end{lemma}
\emph{Proof.}  Set $i=1,2.$ Define $ T^{n} $ the shift operator by $ T^{n}\left(\xi_{i,0},\xi_{i,1},\dots\right)=\left(\xi_{i,n},\xi_{i,n+1},\dots\right) $, for $ n\geq 1 $. Then we get for a fixed $ k\geq 0 $, $$\Pi_{i,n}\left(T^{k}\xi_i\right)=m_{i,k}m_{i,k+1}\cdots m_{i,k+n-1}  .$$
In particular, we have $\Pi_{i,n}=\Pi_{i,n}\left(T^{0}\xi_i\right)$.
Next, we use the method of Grama \emph{et al.}\ \cite{grama2021} to complete the proof.
From  (3.15) in \cite{grama2021}, it is proven that for all $ t>0$ and $n\geq1 $,
\begin{align}\label{v2}
\phi_i(t) &\leq \mathbb{E}\phi_i\left(\frac{t}{\Pi_{i,n}}\right)\prod_{j=0}^{n-1}\Big(p_{1}(\xi_{i,j})+(1-p_{1}(\xi_{i,j}))\beta_{K}\Big) \nonumber\\
 &\quad +\frac{1}{K}\mathbb{E}\left[\phi_{i, T^{n}\xi_i}\left(\frac{t}{\Pi_{i,n}}\right)\left(\mathbb{E}_{T^{n}\xi_i}W_{i, \infty}^{p}\right)\right]\ +\ \mathbb{P}\left(\frac{t}{\Pi_{i,n}}<t_{K}\right),
\end{align}
where $ t_{K}:=(CK)^{-1/(p-1)} $, $C$ and $K$ are positive constants, such that
$$\beta_{K}:=1-(1-1/p)t_{K} \in (0, 1) .$$
From (3.18) in \cite{grama2021}, it  is proven the following inequality
\begin{align}
&\mathbb{E}\left[\phi_{i,T^{n}\xi_i}\left(\frac{t}{\Pi_{i,n}}\right)\mathbb{E}_{T^{n}\xi_i}W_{i,\infty}^{p} \right]\nonumber \\
 \leq\ & \mathbb{E}\phi_i\left(\frac{t}{\Pi_{i,n}}\right)+\sum_{k=0}^{\infty}\mathbb{E}\frac{\phi_i\left(\frac{t}{ \Pi_{i,k+1}(T^{n}\xi_i)
\Pi_{i,n}(\xi_i)}\right)}{\Pi_{i,k}^{p-1 }(T^{n}\xi_i)}\frac{m_{k}^{(p)}\left(T^{n}\xi_i\right)}{m_{k}^{p}\left(T^{n}\xi_i\right)}.\label{v4.3}
\end{align}
Let $p\in(1,2]  $. Given $n, K,$ let  $ \widetilde{N}_{i,n,K} $ be a positive random variable whose distribution is defined by
\begin{align} \mathbb{E}g(\widetilde{N}_{i,n,K})&=\frac{1}{q_{i,n,K}}\left[\mathbb{E}g\left(\frac{1}{\Pi_{i,n}}\right)\prod_{j=0}^{n-1}\Big(p_{1}(\xi_{i,j})+(1-p_{1}(\xi_{i,j}))\beta_{K}\Big)\right. \nonumber\\ &\quad +\left. \frac{1}{K }\mathbb{E}g\left(\frac{1}{\Pi_{i,n}}\right)+\frac{1}{K }\sum_{k=0}^{\infty}\mathbb{E}\frac{g\left( \frac{1}{  \Pi_{i,k+1}\left(T^{n}\xi_i\right)\Pi_{i,n}\left(\xi_i\right)}\right)}{\Pi_{i,k}^{p-1 }\left(T^{n}\xi_i\right)}\frac{m_{i,k}^{(p)}\left(T^{n}\xi_i\right)}{m_{i,k}^{p}\left(T^{n}\xi_i\right)} \right],\label{v4.1}
\end{align}
for any bounded and measurable function $ g $, where  $ q_{i,n,K}$ is the normalizing constant
(to make $ \mathbb{E} g(\widetilde{N}_{i,n,K})=1$, when $ g=1 $) defined by
\begin{align}
	q_{i,n,K}&=\mathbb{E}\bigg[\prod_{j=0}^{n-1}\left(p_{1}(\xi_{i,j})+(1-p_{1}(\xi_{i,j}))\beta_{K}\right)\bigg]   +\frac{1}{K }\left[1+\sum_{k=0}^{\infty}\mathbb{E}\frac{1}{\Pi_{i,k}^{p-1 }\left(T^{n}\xi_i\right)}\frac{m_{i,k}^{(p)}\left(T^{n}\xi_i\right)}{m_{i,k}^{p}\left(T^{n}\xi_i\right)} \right] \nonumber \\
	&=\Big[ \mathbb{E} \left(p_{1}(\xi_{i,0})+(1-p_{1}(\xi_{i,0}))\beta_{K}\right) \Big]^n   +\frac{1}{K }\left[1 +
	\frac{1}{1 - \mathbb{E} m_{i,0}^{-(p-1)}}  \mathbb{E}\Big( \frac{m_{i,0}^{(p)} }{m_{i,0}^{p} }\Big) \right]. \label{defq}
\end{align}
The last line follows by (3.21) in \cite{grama2021}.
Combining \eqref{v2}-\eqref{v4.1} together, we have
\begin{align}\label{v4}
	\phi_i(t)\leq q_{i,n,K}\mathbb{E}\phi_i(\widetilde{N}_{i,n,K}\,t)+\mathbb{P}\left(\frac{t}{\Pi_{i,n}}<t_{K}\right).
\end{align}
Choose $A_i$ such that $\ln A_i>\mu_i$. Clearly, we have $A_i>1$.
When $A_i^{n + 1}\leq t< A_i^{n + 2}$, let $K =K_{i,n}=((n + 1)\ln A_i )^{1+\delta} $,
then we have $\lim_{n\rightarrow \infty} \sqrt[n]{t_{K_{i,n}}} = 1$. By Fuk-Nagaev's inequality (see Corollary 2.5 in \cite{FGL17}
with  $p=2+ \delta,\ V^2=O(n)$ and $C_p=O(n)$), we have for $n$ large enough,
\begin{align}
	\mathbb{P}\Big(\frac{t}{\Pi_{i,n}}<t_{K_{i,n}}\Big)&=\mathbb{P}\Big(\Pi_{i,n}>\frac{t}{t_{K_{i,n}}}\Big)\leq \mathbb{P}\Big(S_{i,n}-n \mu_i > n\big(\ln \frac{A_i}{\sqrt[n]{t_{K_{i,n}}}} -\mu_i \big)\Big) \nonumber \\
&\leq  \exp\Big\{  - c \, n\, c_i^2 \Big\} + C_1\frac{ n\, c_i^{2+\delta} }{n^{2+\delta}}\nonumber \\
&\leq\frac{C_2}{n^{1+\delta}}\label{v4.5},
\end{align}
where $S_{i,n}= \sum_{j=1}^n M_{i, j-1}$ and $ c_i=  \ln \frac{A_i}{\sqrt[n]{t_{K_{i,n}}}} -\mu_i,\ i=1, 2$.
When $t\geq A_i^{2}$, let $n$ be a positive integer such that $  A_i^{n+1}\leq t< A_i^{n+2},$  so we get
$$n>\frac{\ln  t   }{\ln A_i}-2.$$
Thus, by \eqref{v4.5}, for all $  A_i^{n+1}\leq t<  A_i^{n+2},$
\begin{align}
	\mathbb{P}\Big(\frac{t}{\Pi_{i,n}}<t_{K_{i,n}}\Big) &\leq \frac{C}{n^{1+\delta}} < C\left(\frac{\ln  t   }{\ln A_i}-2\right)^{-1-\delta} \nonumber \\
 &\leq \frac{C}{(\ln t)^{1+\delta}}.  \label{v4.6}
\end{align}
Notice that $0\leq  \phi_i(t)\leq 1\ (t>0)$. By \eqref{v4} and \eqref{v4.6}, we have
\begin{equation}\label{v4.128}
	\phi_i(t)\leq\sum_{n=1}^{\infty}q_{i,n,K_{i,n}}  \mathbf{1}_{\{  A_i^{n +1}\leq t<  A_i^{n +2} \} } +\frac{C}{1+(\ln^+ t)^{1+\delta}}   .
\end{equation}
From the definition \eqref{defq} of $q_{i,n,K}$, when $  A_i^{n +1}\leq t<  A_i^{n +2}$, we have
\begin{equation}\label{v4.1df8}
	q_{i,n,K_{i,n}}   \leq   \, \frac{C}{1+(\ln^+ t)^{1+\delta}}  .
\end{equation}
Finally, by \eqref{v4.128} and \eqref{v4.1df8}, we obtain for all $t > 0$,
\begin{equation*}
	\phi_i( t)\leq   \frac{C}{1+(\ln^+ t)^{1+\delta}}.
\end{equation*}
This completes the proof of Lemma \ref{lemma4}.
\hfill\qed

  Grama \emph{et al.}\ \cite{GLE17} (see Theorem 3.1 therein) have established  a bound $ \varphi(t) \leq   C\, t^{-\alpha}, t>0,$
 where $\alpha$ is a positive constant. Their bound is better than the one in
Lemma \ref{lemma4}. However, Theorem 3.1 of Grama \emph{et al.}\ \cite{GLE17}
requires condition  \textbf{A3} which is stronger than   condition  \textbf{A1}.

We have the following result for the $\mathbb{L}^p(\mathbb{P})$ moments of $\ln W_{i, \infty}$ and $\ln W_{i, n}$.
The same result with $ q \in (1, 1+\delta/2)$ has been established by Grama \emph{et al.}\ \cite{GLE17} (cf.\ Lemma 2.3 therein).
The following result  is an improvement on their result by replacing $ q \in (1, 1+\delta/2)$ with $ q \in (1, 1+\delta )$.
 \begin{lemma}\label{lemma4.3}
		Assume that the conditions \textbf{A1} and \textbf{A2} are satisfied. Then for $i=1, 2$ and  $\ q \in (1, 1+\delta)$, the following two inequalities hold
		\begin{equation}\label{eqLiu3}
			\mathbb{E}|\ln W_{i, \infty}|^q<\infty
		\end{equation}
		and
		\begin{equation}\label{eqLiu4}
			\sup_{n\in\mathbb{N}}\mathbb{E}|\ln W_{i,n}|^q <\infty.
		\end{equation}
\end{lemma}
\emph{Proof.}  Set $i=1, 2.$  Consider the following truncation
\begin{equation}\label{eqLiu5}
	\mathbb{E}|\ln W_{i, \infty}|^q=\mathbb{E}|\ln W_{i, \infty}|^q \mathbf{1}_{\{W_{i, \infty}>1\}}+	\mathbb{E}|\ln W_{i, \infty}|^q \mathbf{1}_{\{W_{i, \infty}\leq1\}}.
\end{equation}
For the first term in the right-hand side of the equality above, we have
\begin{equation}\label{eqLiu6}
	\mathbb{E}|\ln W_{i, \infty}|^q \mathbf{1}_{\{W_{i, \infty}>1\}}\leq   C\,  \mathbb{E}W_{i, \infty}<\infty.
\end{equation}
For the second term, we have
\begin{align}
	\mathbb{E}|\ln W_{i, \infty}|^q\mathbf{1}_{\{W_{i, \infty}\leq1\}}&=\, q\int_1^{\infty}\frac{1}{t}(\ln t)^{q-1}\mathbb{P}(W_{i, \infty}\leq t^{-1})\, dt \nonumber \\
	&\leq\, q\, e\int_1^{\infty}\frac{\phi_i(t)}{t}(\ln t)^{q-1}\, dt\nonumber\\
	&=q\, e\left(\int_1^{e}\frac{\phi_i(t)}{t}(\ln t)^{q-1}\, dt+\int_e^{\infty}\frac{\phi_i(t)}{t}(\ln t)^{q-1}\, dt\right). \label{eqLiu7}
\end{align}
The last inequality above can be obtained by Markov's inequality, i.e.,
$$\mathbb{P}(W_{i, \infty}\leq t^{-1})\leq e\,\mathbb{E}e^{-tW_{i, \infty}}=e\,\phi_i(t).$$
Clearly, it holds
\begin{equation}\label{i4.13}
	\int_1^{e}\frac{\phi_i(t)}{t}(\ln t)^{q-1}\, dt<\infty.
\end{equation}
From Lemma \ref{lemma4} and $q<1+\delta$, we have
\begin{equation}\label{eqLiu28}
	\int_e^{\infty}\frac{\phi_i(t)}{t}(\ln t)^{q-1}\, dt\leq\, C \int_e^{\infty}\frac{1}{t(\ln t)^{2+\delta-q}}\, dt<\infty.
\end{equation}
Substituting \eqref{i4.13} and \eqref{eqLiu28} into \eqref{eqLiu7}, we get
\begin{equation}\label{i4.15}
	\mathbb{E}|\ln W_{i, \infty}|^q \mathbf{1}_{\{W_{i, \infty}\leq1\}}<\infty.
\end{equation}
Therefore, by (\ref{eqLiu5}), (\ref{eqLiu6}) and (\ref{i4.15}), we get (\ref{eqLiu3}).

Next, we give a proof for (\ref{eqLiu4}). Since $x\mapsto\left|\ln^q(x) \mathbf{1}_{\{x\leq1\}}\right|, q>1,$ is a non-negative and convex function, by Lemma 2.1 in  \cite{HL12}, we have
$$\sup_{n\in\mathbb{N}}\mathbb{E}\left|\ln W_{i,n}\right|^q\mathbf{1}_{\{W_{i,n}\leq1\}}=\mathbb{E}\left|\ln W_{i, \infty}\right|^q\mathbf{1}_{\{W_{i, \infty}\leq1\}}.$$
With the similar truncation as  $\mathbb{E}\left|\ln W_{i, \infty}\right|^q$, by (\ref{i4.15}), we get
\begin{align*}
		\sup_{n\in \mathbb{N}}\mathbb{E}\left|\ln W_{i,n}\right|^q&=\,\sup_{n\in \mathbb{N}}\left(\mathbb{E}\left|\ln W_{i,n}\right|^q\mathbf{1}_{\{W_{i,n}>1\}}+\mathbb{E}\left|\ln W_{i,n}\right|^q\mathbf{1}_{\{W_{i,n}\leq1\}}\right)\\
		&\leq\,\sup_{n\in \mathbb{N}}\mathbb{E}\left|\ln W_{i,n}\right|^q\mathbf{1}_{\{W_{i,n}>1\}}+\sup_{n\in \mathbb{N}}\mathbb{E}\left|\ln W_{i,n}\right|^q\mathbf{1}_{\{W_{i,n}\leq1\}}\\
		&\leq\,C\,\mathbb{E}W_{i, \infty}+\mathbb{E}\left|\ln W_{i, \infty}\right|^q\mathbf{1}_{\{W_{i, \infty}\leq1\}}<\infty.
\end{align*}
This completes the proof of Lemma \ref{lemma4.3}.
\hfill\qed

The following lemma is a simple consequence of Lemma 2.4 of Grama \emph{et al.}\ \cite{GLE17}.
\begin{lemma}\label{lemma 4.4}
		Assume that the conditions \textbf{A1} and \textbf{A2} are satisfied. Then there exists a constant $\gamma\in(0,1)$, such that
		$$  \mathbb{E}|\ln W_{1,n}-\ln W_{1, \infty}| +   \mathbb{E}|\ln W_{2,m}-\ln W_{2, \infty}|\leq C\, \gamma^{m \wedge n}.$$		
\end{lemma}

In the proof of  Theorem \ref{th02}, the following lemma plays an important role.
 \begin{lemma}\label{lemma4.5}
 		Assume that the conditions \textbf{A1} and \textbf{A2} are satisfied.  Let $\delta'$ be a constant such that $\delta' \in (0, \delta).$ Then  for all $ x \in \mathbb{R} $,
 		\begin{equation}\label{7}
 			\mathbb{P}\bigg(R_{m,n} \leq
 			x, \sum_{i=1}^{n+m}  \eta_{ m,n, i}\geq x\bigg) \leq  \frac{C}{(m \wedge n)^{\delta / 2}}\frac{1}{1+|x|^{1+\delta'\ }}
 		\end{equation}
 		and
 		\begin{equation}\label{11}
 			\mathbb{P}\bigg(R_{m,n} \geq x, \sum_{i=1}^{n+m}  \eta_{ m,n, i}\leq x\bigg) \leq
 			\frac{C}{(m \wedge n)^{\delta / 2}}\frac{1}{1+|x|^{1+\delta' \ }}.
 		\end{equation}		
 \end{lemma}
\emph{Proof.} We only give a proof for inequality (\ref{7}), as inequality  \eqref{11} can be proved in the same way. Without loss of generality, we may assume that $  m\leq n$.

First, we show that inequality (\ref{7}) holds  for $x  \leq - C m^{1/2}$ with some positive constant $C$.
Recall
\begin{equation} \nonumber
R_{m,n} =  \frac{  \ln  Z_{1,n}  -n \mu_1   }{ n V_{m,n, \rho}}  -   \frac{   \ln  Z_{2,m}  - m\mu_2 }{m V_{m,n, \rho}},\ \ \quad m,n \in \mathbb{N}.
\end{equation}
Then  we have for all $x  \in \mathbb{R}$,
\begin{equation*}
		\mathbb{P}\bigg(R_{m,n} \leq x, \sum_{i=1}^{n+m}  \eta_{ m,n, i}\geq x\bigg)\leq\mathbb{P}\bigg(R_{m,n} \leq x\bigg)\leq P_1 + P_2,
\end{equation*}
where
$$ P_1 = \mathbb{P}\bigg(   \frac{  \ln  Z_{1,n}  -n \mu_1   }{ n V_{m,n, \rho}}  \leq  \frac x 2   \bigg) \quad \textrm{and}\quad
  P_2  = \mathbb{P}\bigg(  -   \frac{   \ln  Z_{2,m}  - m\mu_2 }{m V_{m,n, \rho}} \leq  \frac x 2   \bigg) .$$
Since $Z_{1,n} \geq 1$ $\mathbb{P}$-a.s.\ and $ V_{m,n, \rho}  \asymp m^{-1/2}$ as $m\rightarrow \infty$, there exists a positive constant $C$ such that    $$ \frac{  \ln  Z_{1,n}  -n \mu_1   }{ n V_{m,n, \rho}}  > - \frac{   \mu_1   }{ \ \ V_{m,n, \rho}} > - \frac12 C  m^{1/2}\ \ \  \mathbb{P}\textrm{-a.s.},  $$
and  thus $P_1=0$ for all $x \leq -C m^{1/2}.$  For $P_2,$ by Lemma \ref{lemma4.1}, Markov's inequality and the fact $ \mathbb{E} W_{2,m} =1$, we have for all $x \leq -C m^{1/2},$
\begin{align*}
	P_2  &\leq   \mathbb{P}\bigg(    \sum_{j=1}^{m} \eta_{m,n,n+j}  \leq -\frac{ |x|} 4  \bigg) + \mathbb{P}\bigg(    \frac{\ln W_{2,m}}{ m\, V_{m,n,\rho} \ }  \geq \frac{ |x|} 4  \bigg) \\
&\leq   \mathbb{P}\bigg(    \sum_{j=1}^{m} \eta_{m,n,n+j}  \leq - \frac{ |x|} 4  \bigg) +  \exp\bigg\{ - \frac{ |x|} 4 m\, V_{m,n,\rho}  \bigg\} \mathbb{E} W_{2,m} \\
&\leq   \frac{C}{m^{\delta / 2}}\frac{1}{1+|x|^{2+\delta}}.
\end{align*}
Hence, inequality  (\ref{7})  holds for all $x  \leq - C m^{1/2}. $

Next, we show that inequality (\ref{7}) holds  for all $x\geq C m^{1/2}$.
For all $x\geq 0,$ we have
\begin{equation*}
		\mathbb{P}\bigg(R_{m,n} \leq
 x, \sum_{i=1}^{n+m}  \eta_{ m,n, i}\geq x\bigg)  \leq  \mathbb{P}\bigg(  \sum_{i=1}^{n+m}  \eta_{ m,n, i}\geq x\bigg)  .
\end{equation*}
Applying Lemma \ref{lemma4.1} to the right-hand side of the last inequality,
we have  for all $x\geq  0,$
\begin{align*}
	&\mathbb{P}\bigg(R_{m,n} \leq x, \sum_{i=1}^{n+m}  \eta_{ m,n, i}\geq x\bigg)\\
	\leq\ &  1-\Phi (x) + \frac{C_1}{ 1+|x| ^{2+\delta}} \Bigg( \sum_{i=1}^{ m} \mathbb{E}|\eta_{ m,n, i} + \eta_{ m,n, n+i} |^{2+\delta}  +\sum_{i=m+1}^{n} \mathbb{E}|\eta_{ m,n, i}  |^{2+\delta}  \Bigg) .
\end{align*}
Using the inequality $$ (a + b)^{2+\delta} \leq 2^{1+\delta} (|a| ^{2+\delta}  +  | b| ^{2+\delta} ), \ \ \ \ \  \ \  a, b \in \mathbb{R},$$
we deduce that
\begin{align}
&\sum_{i=1}^{ m} \mathbb{E}|\eta_{ m,n, i} + \eta_{ m,n, n+i} |^{2+\delta}  +\sum_{i=m+1}^{n} \mathbb{E}|\eta_{ m,n, i}  |^{2+\delta}\nonumber\\
  & \ \ \ \ \ \ \leq\ C_1   \Bigg( \sum_{i=1}^{n} \mathbb{E}|\eta_{ m,n, i} |^{2+\delta}  +\sum_{i= 1}^{m } \mathbb{E}|\eta_{ m,n, n+i}  |^{2+\delta}  \Bigg)	.\label{ines32}
\end{align}
Hence, we get for all $x\geq 0,$
\begin{align*}
	&\mathbb{P}\bigg(R_{m,n} \leq
 x, \sum_{i=1}^{n+m}  \eta_{ m,n, i}\geq x\bigg) \\
 \leq\ & 1-\Phi (x) +  \ \frac{C_1}{ 1+|x| ^{2+\delta}} \Bigg( \sum_{i=1}^{n} \mathbb{E}|\eta_{ m,n, i} |^{2+\delta}  +\sum_{i= 1}^{m } \mathbb{E}|\eta_{ m,n, n+i}  |^{2+\delta}  \Bigg).
\end{align*}
Notice that $V_{m,n,\rho} \asymp  m^{-1/2}$ as $m  \rightarrow \infty$.
Then,  by the  inequalities
\begin{equation}\label{norb}
\frac 1{\sqrt{2 \pi} ( 1+x)  }e^{- x^2/2}\leq
1-\Phi \left( x\right)   \leq  \frac 1{\sqrt{ \pi} ( 1+x)  }e^{- x^2/2}, \ \   x \geq 0,
\end{equation}
 we have  for all $x\geq C m^{1/2},$
\begin{align*}
 	\mathbb{P}\bigg(R_{m,n} \leq x, \sum_{i=1}^{n+m}  \eta_{ m,n, i}\geq x\bigg)
 &\leq\frac{C_2}{m^{\delta / 2}}\frac{1}{1+|x|^{2+\delta}} + \frac{C_3}{ 1+|x| ^{2+\delta}} \Bigg( \frac{n}{n^{2+\delta} m^{-1 -\delta/2} } + \frac{m}{m^{2+\delta} m^{-1 -\delta/2} }    \Bigg)\\
 &\leq\frac{C_4}{m^{\delta / 2}}\frac{1}{1+|x|^{2+\delta}} .
\end{align*}
Thus,  inequality (\ref{7}) holds  for all $x\geq C m^{1/2} .$

To end the proof of lemma,  we only need to show that (\ref{7}) holds for all $  |x| < C m^{1/2}$.   Consider the following notations
for all $0\leq k \leq m-1$,
\begin{align}
 Y_{m,n,k}&=\sum_{i=k+1}^{n} \eta_{ m,n, i}+\sum_{j=k+1}^{m} \eta_{ m,n, n+j}, \label{Y1}\\
\widetilde{Y}_{m,n, k}&=Y_{m,n,0}- Y_{m,n,k} ,\label{Y2}\\
H_{m,n,k}&=\frac{\ln W_{1, k}}{ n\, V_{m,n,\rho}  \ } - \frac{\ln W_{2,k}}{ m\, V_{m,n,\rho}}, \label{Y3} \\
 D_{m, n , k}&=\frac{\ln W_{1, n}}{ n\, V_{m,n,\rho}  \ } - \frac{\ln W_{2,m}}{ m\, V_{m,n,\rho} \ }-H_{m,n,k}. \label{Y4}
\end{align}
Set $ \alpha_{m}= m^{-\delta/2}$ and $k=[ m^{1-\delta/2} \,]$, where   $ [t]$ stands for the largest integer less than $t$.
From \eqref{4.8}, we deduce that for all $  x \in \mathbb{R}$,
\begin{align}
&\mathbb{P}\bigg(R_{m,n} \leq x, \sum_{i=1}^{n+m}  \eta_{ m,n, i}\geq x\bigg)\nonumber\\
=\ & \mathbb{P}\bigg(Y_{m,n,0}+H_{m,n,k}+D_{m, n , k}\leq x, Y_{m,n,0}\geq x\bigg)  \nonumber  \\
\leq\ &\mathbb{P}\bigg(Y_{m,n,0}+H_{m,n,k}\leq x+ \alpha_{m}, Y_{m,n,0}\geq x\bigg)+ \mathbb{P}\bigg(|D_{m, n , k}|\geq \alpha_{m}\bigg). \label{4.38}
\end{align}
For the tail probability $\mathbb{P}\left(|D_{m, n , k}|\geq \alpha_{m}\right)$, by Markov's inequality and
Lemma \ref{lemma 4.4}, there exists a constant $ \gamma \in (0,1)$ such that for all  $ - m < x < m$,
\begin{align}\label{4.39}
	\mathbb{P}\left(\left|D_{m, n, k}\right|>\alpha_{m}\right)
	& \leq \frac{\mathbb{E}\left|D_{m, n,k}\right|}{\alpha_{m}}\nonumber\\
	&=\frac{m^{ \delta/2} }{V_{m,n,\rho}}\mathbb{E}\Bigg|
	\left[\left(\frac{\ln W_{1, n}}{n}-\frac{\ln W_{1, \infty}}{n}\right) - \left(\frac{\ln W_{2,m}}{m}-\frac{\ln W_{2,\infty}}{m}\right)\right]\nonumber\\
	&\quad -\left[\left(\frac{\ln W_{1, k}}{n}-\frac{\ln W_{1, \infty}}{n}\right) - \left(\frac{\ln W_{2,k}}{m}-\frac{\ln W_{2, \infty}}{m}\right)\right]\Bigg|  \nonumber\\
&\leq\frac{m^{ \delta/2} }{V_{m,n,\rho}}\Bigg(
	\mathbb{E} \left|\frac{\ln W_{1, n}}{n}-\frac{\ln W_{1, \infty}}{n}\right| + \mathbb{E}  \left|\frac{\ln W_{2,m}}{m}-\frac{\ln W_{2,\infty}}{m}\right| \nonumber\\
	&\quad +  \mathbb{E}  \left|\frac{\ln W_{1, k}}{n}-\frac{\ln W_{1, \infty}}{n}\right|  + \mathbb{E} \left|\frac{\ln W_{2,k}}{m}-\frac{\ln W_{2, \infty}}{m}\right| \Bigg)  \nonumber\\
	&\leq C_1 \, m^{(1+  \delta)/2 } \bigg( \frac1n \gamma ^{n}   + \frac1m \gamma ^{m}+ \frac1n \gamma ^{k} + \frac1m \gamma ^{k} \bigg)  \nonumber  \\
& \leq \frac{C_2 }{m^{\delta / 2}} \frac{1}{1+|x|^{2+\delta}}.
\end{align}
Next, we give an estimation for  the first term at the right-hand side of \eqref{4.38}. Let
\begin{equation}\label{inds25a}
  G_{m, n,k}(x)=\mathbb{P}\left(Y_{m, n,k} \leq x\right)\ \text{ and } \ v_{k}(d s, d t)=\mathbb{P}\left(\widetilde{Y}_{m,n, k} \in d s, H_{m,n,k} \in d t\right).
\end{equation}
Since $ Y_{m, n,k} $ and $(\widetilde{Y}_{m,n, k},H_{m,n,k})  $ are independent, we have
\begin{align}
	&\mathbb{P}\Big(Y_{m,n,0}+H_{m,n,k} \leq x+\alpha_{m},  Y_{m,n,0} \geq x\Big)\nonumber \\
	=\ & \mathbb{P}\left(Y_{m,n,k}+\widetilde{Y}_{m,n, k}  +H_{m,n,k} \leq x+\alpha_{m},  Y_{m,n,k}+ \widetilde{Y}_{m,n, k}  \geq x\right) \nonumber\\
	=\ & \iint \mathbb{P}\left(Y_{m,  n,k}+s+t \leq x+\alpha_{m},  Y_{m,n,k}+s \geq x\right) v_{k}(ds,  dt) \nonumber\\
	=\ &  \iint \mathbf{1}_{\{t \leq \alpha_{m}\}}\Big(G_{m,  n,k}\left(x-s-t+\alpha_{m}\right)-G_{m,  n,k}(x-s)\Big) v_{k}(ds,  dt).\label{4.40}
\end{align}
Denote $  C_{m,n,k}^2=  \textrm{Var} (Y_{m,n,k}),  $ then it holds $ C_{m,n,k}=  1 + O(k/n) \nearrow 1  $ as $ m \rightarrow \infty.$
By Lemma \ref{lemma4.1} and (\ref{ines32}),  we obtain for all $x \in \mathbb{R},$
\begin{align*}
&\left|\mathbb{P}\left(  \frac{ Y_{m,n,k} }{C_{m,n,k} }  \leq \frac{x}{C_{m,n,k}}\right)-\Phi\left(\frac{x}{C_{m,n,k}}\right)\right|\\
\leq \ & \frac{C_1}{1+|x/C_{m,n,k}|^{2+\delta}}\left(\sum_{j=k+1}^{m}\mathbb{E}\left|\frac{\eta_{ m,n, i}}{ C_{m,n,k}} + \frac{\eta_{ m,n, n+ i}}{ C_{m,n,k}}\right|^{2+\delta} + \sum_{i=m+1}^{n}\mathbb{E}\left|\frac{\eta_{ m,n, i}}{ C_{m,n,k}}\right|^{2+\delta}
\right)\\
\leq\ & \frac{C_2}{1+|x |^{2+\delta}}\left( \sum_{j=k+1}^{m}\mathbb{E}\left|\frac{M_{2,j-1}-\mu_2}{mV_{m,n,\rho} }\right|^{2+\delta} + \sum_{i=k+1}^{n}\mathbb{E}\left|\frac{M_{1,i-1}-\mu_1}{nV_{m,n,\rho} }\right|^{2+\delta}\right)\\
\leq\ & \frac{C_3}{1+|x|^{2+\delta}} \Bigg(\frac{m}{m^{2+\delta} m^{-1 -\delta/2} }+ \frac{n}{n^{2+\delta} m^{-1 -\delta/2} }   \Bigg) \\
\leq\ & \frac{C_4}{m^{\delta / 2}}\frac{1}{1+|x|^{2+\delta}} .
\end{align*}
By the last inequality, we deduce that for all  $x \in \mathbb{R},$
\begin{align*}
\left|G_{m,  n,k}(x)-\Phi(x)\right|
&\leq \left|\mathbb{P}\left( \frac{ Y_{m,n,k} }{C_{m,n,k} }\leq \frac{x}{C_{m,n,k}}\right)-\Phi\left(\frac{x}{C_{m,n,k}}\right)\right|
  +  \left|\Phi\left(\frac{x}{C_{m,n,k}}\right)-\Phi\left(x\right)\right|\nonumber\\
&\leq  \frac{C}{m^{\delta / 2}} \frac{1}{1+|x|^{2+\delta}}+\exp\left\{-\frac{x^{2}}{2}\right\}\left|\frac{x}{C_{m,n,k}}-x\right|\nonumber\\
&\leq    \frac{C}{m^{\delta / 2}} \frac{1}{1+|x|^{2+\delta}} +  C \frac{k}{n} \frac{1}{1+|x|^{2+\delta}}  \nonumber   \\
&\leq    \frac{C}{m^{\delta / 2}} \frac{1}{1+|x|^{2+\delta}}. \nonumber
\end{align*}
Therefore,   we have for all $x \in \mathbb{R},$
\begin{equation}\label{4.43}
  \mathbb{P}\Big(Y_{m,n,0}+H_{m,n,k} \leq x+\alpha_{m},  Y_{m,n,0} \geq x\Big) \leq  J_{1}+J_{2} +J_3,
\end{equation}
where
$$J_{1}=\iint \mathbf{1}_{\{t \leq \alpha_{m}\}}\left|\Phi\left(x-s-t+\alpha_{m}\right)-\Phi(x-s)\right| v_{k}(d s,  d t),$$
$$ J_{2}=\frac{C}{m^{\delta / 2}} \iint  \mathbf{1}_{\{t \leq \alpha_{m}\}}  \frac{1}{1+|x-s|^{2+\delta}} v_{k}(d s,  d t) $$
and
$$ J_{3}=\frac{C}{m^{\delta / 2}}\iint \mathbf{1}_{\{t \leq \alpha_{m}\}} \frac{1}{1+|x-s-t|^{2+\delta}} v_{k}(d s,  d t). $$
For $J_{1}$, by the mean value theorem, we have for all $x \in \mathbb{R},$
\begin{align*}
&\mathbf{1}_{\{t \leq \alpha_{m}\}}\left|\Phi\left(x-s-t+\alpha_{m}\right)-\Phi(x-s)\right| \\
 \leq \ & C|\alpha_{m}-t |\exp\left\{-\frac{x^{2}}{8}\right\} + \,  |\alpha_{m}-t | \mathbf{1}_{\{|s| \geq 1+  \frac{1}{4 }|x| \}} +   \,  |\alpha_{m}-t | \mathbf{1}_{\{|t| \geq 1+  \frac{1}{4 }|x| \}}
\end{align*}
and
\begin{equation}\label{4.44}
	J_{1}\leq J_{11}+ J_{12}+ J_{13},
\end{equation}
where
\begin{align*}
	J_{11}&= C\iint |\alpha_{m}-t |\exp\left\{-\frac{x^{2}}{8}\right\} v_{k}(d s,  d t),\\
J_{12}&=  \iint |\alpha_{m}-t |\mathbf{1}_{\{|s| \geq 1+  \frac{1}{4 }|x| \}}  v_{k}(d s,  d t)
\end{align*}
and
\begin{equation*}
	J_{13}= \iint |\alpha_{m}-t |\mathbf{1}_{\{|t| \geq 1+  \frac{1}{4 }|x| \}} v_{k}(d s,  d t).
\end{equation*}
By Lemma \ref{lemma4.3}, it is obvious that  for all $x \in \mathbb{R},$
\begin{align*}
	J_{11}& \leq   C_1 \exp\left\{-\frac{x^{2}}{8}\right\}  \bigg(\alpha_{m} + \mathbb{E}| H_{m,n,k}| \bigg)\\
&\leq   \frac{C_2}{m^{\delta/ 2}} \frac{1}{1+|x|^{2+\delta}}.
\end{align*}
For $J_{12},$ we have the following estimation  for all $x \in \mathbb{R},$
\begin{equation*}
	J_{12} \leq \alpha_{m} \mathbb{P}\bigg( |\widetilde{Y}_{m,n, k} | \geq 1+  \frac{1}{4 }|x| \bigg)   + \mathbb{E}| H_{m,n,k}|\mathbf{1}_{\{|\widetilde{Y}_{m,n, k} | \geq 1+  \frac{1}{4 }|x| \}} . \nonumber
\end{equation*}
Denote $\widetilde{C}_{m,n,k}^2= \textrm{Var}(\widetilde{Y}_{m,n, k} )$, then we have $ \widetilde{C}_{m,n,k}^2\asymp \frac{1}{m^{\delta/2}}.$
Let $\delta' \in (0, \delta).$
By Lemma \ref{lemma4.1}, we deduce that for all $x \in \mathbb{R},$
\begin{align}
 \mathbb{P}\bigg( |\widetilde{Y}_{m,n, k} | \geq 1+  \frac{1}{4 }|x| \bigg)
&=\mathbb{P}\left(  \frac{\widetilde{Y}_{m,n, k} }{\widetilde{C}_{m,n,k} } \geq\frac{1+|x|/4}{ \widetilde{C}_{m,n,k} } \right)+\mathbb{P}\left(   \frac{\widetilde{Y}_{m,n, k} }{\widetilde{C}_{m,n,k} }  \leq -\frac{1+|x|/4}{ \widetilde{C}_{m,n,k} } \right)\nonumber\\
&\leq 1-\Phi\left(\frac{1+|x|/4}{ \widetilde{C}_{m,n,k} }\right)+\Phi\left(-\frac{1+|x|/4}{ \widetilde{C}_{m,n,k} }\right)\nonumber\\
&\quad + \frac{C}{1+ \Big| \frac{1+|x|/4}{ \widetilde{C}_{m,n,k} }\Big|^{2+\delta}} \sum_{i=1}^{k}\mathbb{E} \Big| \frac{\eta_{ m,n, i}+  \eta_{ m,n, n+i}}{\widetilde{C}_{m,n,k}} \Big|^{2+\delta} \nonumber\\
&\leq  \frac{C_1}{1+|x|^{2+\delta}}    \frac{1 }{m^{\delta(2+\delta)/4}} \frac{1 }{k^{\delta/2}}      \nonumber \\
&\leq  \frac{C_2}{1+|x|^{2+\delta}}   \frac{1 }{m^{\delta }}  \label{4.46}
\end{align}
and, by H\"{o}lder's inequality with $\tau= 1+ \frac{ \delta + \delta'}{2+2\delta-\delta'} $ and  $\iota $   satisfying  $\frac{1}{\tau} + \frac{1}{\iota}=1$, it holds for all $|x| \leq C m^{1/2},$
\begin{align*}
 \mathbb{E}| H_{m,n,k}|\mathbf{1}_{\{|\widetilde{Y}_{m,n, k} | \geq 1+  \frac{1}{4 }|x| \}}  & \leq    \Big(\mathbb{E}| H_{m,n,k}|^ \tau \Big)^{1/\tau}   \mathbb{P}\bigg( |\widetilde{Y}_{m,n, k} | \geq 1+  \frac{1}{4 }|x| \bigg)  ^{1/\iota} \\
 &\leq C \frac{1}{m^{1/2}}  \bigg( \frac{C_1}{1+|x|^{2+\delta}}   \frac{1 }{m^{\delta }}  \bigg)  ^{1/\iota}\\
 &\leq  \frac{C }{m^{\delta / 2}} \frac{1}{1+|x|^{1+\delta' \ }}.
\end{align*}
Hence, we have for all $|x| \leq C m^{1/2},$
\begin{equation*}
	J_{12}\leq \frac{C_3}{m^{\delta / 2}} \frac{1}{1+|x|^{1+\delta'}}.
\end{equation*}
For $J_{13},$  we have for all $x \in \mathbb{R}$,
\begin{equation*}
	J_{13} \leq  \alpha_{m}  \mathbb{P}\bigg( |H_{m,n,k}| \geq 1+  \frac{1}{4 }|x| \bigg) +   \mathbb{E}| H_{m,n,k}|\mathbf{1}_{\{|H_{m,n,k} | \geq 1+  \frac{1}{4 }|x| \}}  .
\end{equation*}
By Lemma \ref{lemma4.3} with $p'=1 + \delta/2$  and  Markov's inequality, we deduce that  for all $|x| \leq C m^{1/2},$
\begin{align}
 \mathbb{P}\bigg( |H_{m,n,k} | \geq 1+  \frac{1}{4 }|x| \bigg) & \leq   \frac{4^{p'}}{1+ |x|^{p'}}  \mathbb{E} |H_{m,n,k} |^{p'}  \ \leq \   \frac{C}{1+ |x|^{p'}}  \frac{1}{ m^{p'/2}}  \nonumber \\
  & \leq   \frac{C}{1+ |x|^{2+\delta}} ,  \label{ines48}
\end{align}
and, by Lemma \ref{lemma4.3} with $p'' =\frac{1}{2}(\delta+\delta') ,$
\begin{align*}
 \mathbb{E}| H_{m,n,k}|\mathbf{1}_{\{|H_{m,n,k} |
  \geq  1+  \frac{1}{4 }|x| \}} \leq \frac{C_1}{1+ |x|^{p''}}    \mathbb{E}| H_{m,n,k}|^{1+p''}     &\leq    \frac{C_2 }{1+ |x|^{p''}}  \frac{1}{ m^{(1+p'')/2}}\\
    &\leq  \frac{C_3}{m^{\delta / 2}} \frac{1}{1+|x|^{1+\delta'}}  .
\end{align*}
Hence, we have  for all $|x| \leq C m^{1/2}, $
\begin{equation*}
	J_{13} \leq     \frac{C}{m^{\delta / 2}} \frac{1}{1+|x|^{1+\delta'}}  .
\end{equation*}
Returning to (\ref{4.44}), we get  for all $|x| \leq C m^{1/2}, $
\begin{equation} \label{fgsg43}
	J_{1} \leq  \frac{C}{m^{\delta / 2}} \frac{1}{1+|x|^{1+\delta'}\ } .
\end{equation}
Next, we consider $J_{2}$. By an argument similar to the proof of \eqref{4.46}, we have for all $x \in \mathbb{R},$
\begin{align}
J_{2}&=\frac{C_1}{m^{\delta / 2}}\iint  \mathbf{1}_{\{t \leq \alpha_{m}\}}  \frac{1}{1+|x-s|^{2+\delta}} v_{k}(d s,  d t) \nonumber  \\
  &\leq\frac{C_1}{m^{\delta / 2}}\left(\int_{|s|< 1+ |x|/2}\frac{1}{1+|x-s|^{2+\delta}} v_{k}(d s)  + \int_{|s| \geq 1    + |x|/2}\frac{1}{1+|x-s|^{2+\delta}} v_{k}(d s) \right)  \nonumber \\
	& \leq  \frac{C_2}{m^{\delta / 2}}\left[  \frac{1}{1+|x/2|^{2+\delta}}+\mathbb{P}\bigg( \Big| \frac{\widetilde{Y}_{m,n,k} }{\widetilde{C}_{m,n,k} } \Big|  > \frac{1+|x|/2}{ \widetilde{C}_{m,n,k} } \bigg)\right]\nonumber  \\
	&\leq  \frac{C_3}{m^{\delta / 2}}\frac{1}{1+|x|^{2+\delta}}. \label{fgsg44}
\end{align}
For $J_{3}$, by some arguments similar to that of \eqref{4.46} and \eqref{ines48}, we have  for all $|x| \leq C m^{1/2}, $
\begin{align}
J_{3}&=  \frac{C_1}{m^{\delta / 2}}\iint \mathbf{1}_{\{t \leq \alpha_{m}\}} \frac{1}{1+|x-s-t|^{2+\delta}} v_{k}(d s,  d t)\nonumber  \\
  &\leq\frac{C_1}{m^{\delta / 2}}\left(\iint_{|s+t|\leq 2+ |x|/2}\frac{1}{1+|x/2|^{2+\delta}} v_{k}(d s,  d t)  \right. \nonumber \\
  & \left. \quad+\iint_{|s| > 1+ |x|/4} v_{k}(d s,  d t)+ \iint_{|t| > 1+ |x|/4} v_{k}(d s,  d t) \right)  \nonumber \\
	&\leq  \frac{C_2}{m^{\delta / 2}}\left[  \frac{1}{1+|x/2|^{2+\delta}}+\mathbb{P}\bigg(  \Big| \frac{\widetilde{Y}_{m,n,k} }{\widetilde{C}_{m,n,k} } \Big|  > \frac{1+|x|/4}{ \widetilde{C}_{m,n,k} } \bigg) + \mathbb{P}\bigg(   |H_{m,n,k}|  > 1+ \frac{|x|}{4}    \bigg) \right]\nonumber  \\
	&\leq  \frac{C_3}{m^{\delta / 2}}\frac{1}{1+|x|^{2+\delta}}. \label{fgsg45}
\end{align}
Applying  the inequalities \eqref{fgsg43}-\eqref{fgsg45}  to \eqref{4.43}, we get  for all $|x| \leq C m^{1/2}, $
\begin{equation}\label{u4.40}
	\mathbb{P}\Big(Y_{m,n,0}+H_{m,n,k}\leq x+ \alpha_{m}, Y_{m,n,0}\geq x\Big) \leq    \frac{C}{m^{\delta / 2}}\frac{1}{1+|x|^{1+\delta'}}.
\end{equation}
Combining  (\ref{4.38}),  (\ref{4.39}) and \eqref{u4.40} together, we get  (\ref{7})  for all $|x| \leq C m^{1/2}. $
This completes the proof of Lemma \ref{lemma4.5}.
\hfill\qed

\subsection{Proof of Theorem \ref{th02}}
We are now in a position to end the proof of Theorem \ref{th02}.
By  Lemma \ref{lemma4.1} and  the fact $V_{m,n, \rho}  \asymp  \sqrt{m^{-1 } + n^{-1 } } $, we have for all $x \in \mathbb{R},$
\begin{align}
&\left|\mathbb{P}\bigg(\sum_{i=1}^{n+m}  \eta_{ m,n, i}  \leq x\bigg)-\Phi(x)\right| \nonumber \\
\leq&\  \frac{C_1}{ 1+|x| ^{2+\delta}}\Bigg( \sum_{i=1}^{ m} \mathbb{E}|\eta_{ m,n, i} + \eta_{ m,n, n+i} |^{2+\delta}  +\sum_{i=m+1}^{n} \mathbb{E}|\eta_{ m,n, i}  |^{2+\delta}  \Bigg) \nonumber \\
	\leq&\   \frac{C_2}{ 1+|x| ^{2+\delta}} \sum_{i=1}^{m+n} \mathbb{E}\left| \eta_{ m,n, i}\right|^{2+\delta} \nonumber  \\
 \leq&\ \frac{C_3}{ 1+|x| ^{2+\delta}}  \Bigg( \frac{n}{n^{2+\delta} \, ( \frac1n    +  \frac1m )^{(2+\delta)/2}   }  + \frac{m}{m^{2+\delta} \, ( \frac1n    +  \frac1m )^{(2+\delta)/2}   }   \Bigg) \nonumber \\
 \leq&\ \frac{C}{ (m \wedge n)^{ \delta/2}    } \frac{1}{ 1+|x| ^{2+\delta}} .  \label{inds35}
\end{align}
Notice that
\begin{align*}
	\mathbb{P}\left(R_{m,n} \leq x\right)
	&= \mathbb{P}\Big( R_{m,n} \leq x, \sum_{i=1}^{n+m}  \eta_{ m,n, i}  \leq x\Big)+\mathbb{P}\Big(R_{m,n} \leq x, \sum_{i=1}^{n+m}  \eta_{ m,n, i}  >x\Big) \nonumber\\
	 &=  \mathbb{P}\Big( \sum_{i=1}^{n+m}  \eta_{ m,n, i}  \leq x\Big) -\mathbb{P}\Big(R_{m,n} >x, \sum_{i=1}^{n+m}  \eta_{ m,n, i} \leq x\Big)\nonumber\\
	&\quad +\mathbb{P}\Big(R_{m,n} \leq x, \sum_{i=1}^{n+m}  \eta_{ m,n, i}  >x\Big).
\end{align*}
Applying (\ref{inds35}) to the last equality, we deduce that for  all $x \in \mathbb{R},$
\begin{align*}
	\Big|\mathbb{P}\left(R_{m,n} \leq x\right) -\Phi(x)\Big|
 &\leq \frac{C}{ (m \wedge n)^{ \delta/2}    } \frac{1}{ 1+|x| ^{2+\delta}}+ \mathbb{P}\Big(R_{m,n} >x, \sum_{i=1}^{n+m}  \eta_{ m,n, i} \leq x\Big)\\
	&\quad + \mathbb{P}\Big(R_{m,n} \leq x, \sum_{i=1}^{n+m}  \eta_{ m,n, i}  >x\Big) .
\end{align*}
By Lemma \ref{lemma4.5}, it follows that for  all $x \in \mathbb{R},$
\begin{equation*}
	\Big|\mathbb{P}\left(R_{m,n} \leq x\right) -\Phi(x)\Big|
  \leq  \frac{C}{ (m \wedge n)^{ \delta/2}    } \frac{1}{ 1+|x| ^{1+\delta'}}.
\end{equation*}
This completes the proof of Theorem \ref{th02}.
\hfill\qed

\subsection{Preliminary Lemmas for   Theorem \ref{th2.2} }   \label{sec4}

 To prove   Theorem \ref{th2.2}, we shall make use of the following lemma (see Theorem 3.1 of Grama \emph{et al.}\ \cite{GLE17}).
The  lemma
shows that the conditions \textbf{A3} and \textbf{A4} imply the existence of a harmonic moment  of positive order $\alpha >0.$
\begin{lemma}\label{lemma1}
Assume that the conditions \textbf{A3} and \textbf{A4} are satisfied.   There exists a constants $ a_0>0$   such that for all  $\alpha \in (0, a_0),$
the following inequalities hold
 \begin{equation}\label{eqLidsu3}
   \mathbb{E}W_{1, \infty}^{-\alpha} + \mathbb{E}W_{2, \infty}^{-\alpha}< \infty
\end{equation}
and
\begin{equation}\label{eqLisu4}
		\sup_{n\in\mathbb{N}}  \big( \mathbb{E}W_{1, n}^{-\alpha } + \mathbb{E}W_{2, n}^{-\alpha } \big) <\infty.
	\end{equation}
\end{lemma}
\emph{Proof.}    We give an alternative proof for Theorem 3.1 of Grama \emph{et al.}\ \cite{GLE17}.
Let $i=1, 2.$
By the fact that $$W_{i, \infty}^{-\alpha}=\frac{1}{\Gamma\left(\alpha\right)} \int_{0}^{\infty} e^{-tW_{i, \infty} } t^{\alpha-1} d t,$$
we have
\begin{align}
	\mathbb{E}W_{i, \infty}^{-\alpha}
	&=\frac{1}{\Gamma\left(\alpha\right)} \int_{0}^{\infty} \phi_i(t) t^{\alpha-1} dt\nonumber\\
	&=\frac{1}{\Gamma\left(\alpha\right)}\left( \int_{0}^{1} \phi_i(t) t^{\alpha-1} dt+ \int_{1}^{\infty} \phi_i(t) t^{\alpha-1} dt\right), \label{i4.45}
\end{align}
where $\Gamma$ is the gamma function. For the first term in the above bracket, since $0\leq \phi_i(t)\leq 1 \ (t \geq 0)$, then  we have for all $a_{0}>0$,
\begin{equation}\label{i4.46}
	\int_{0}^{1} \phi_i(t) t^{a_{0}-1} dt\leq C \int_{0}^{1}  t^{a_{0}-1} dt<\infty.
\end{equation}
For the second term, it suffices  that there exists a positive constant $a_0$ such that for all $t>0 $,
\begin{equation}\label{u4.47}
	\phi_i(t)\leq  \frac{C}{1+t^{a_0}}.
\end{equation}
To this end, we use  the method of Grama \emph{et al.}\,\cite{grama2021}.
Let $\widetilde{N}_{i,n,K}$ be defined as in (\ref{v4.1}).
By \eqref{v4}, we have for all $t>0$,
\begin{equation}\label{vs4}
	\phi_i(t)\leq q_{i,n,K}\mathbb{E}\phi_i(\widetilde{N}_{i,n, K}\,t)+\mathbb{P}\Big(\frac{t}{\Pi_{i,n}}<t_{K}\Big).
\end{equation}
Let $p\in(1,2]$. By \eqref{v4.1}, we have
\begin{align*}
	q_{i,n,K}\mathbb{E}\widetilde{N}_{i,n,K}^{-a}
	&=\mathbb{E}\left[\Pi_{i,n}^{a}\prod_{j=0}^{n-1}\left(p_{1}(\xi_{i,j})+(1-p_{1}(\xi_{i,j}))\beta_{K}\right)\right]+\frac{1}{K }\,\mathbb{E}\Pi_{i,n}^{a} \\  &\quad +\frac{1}{K }\mathbb{E}\left[\Pi_{i,n}^{a}(\xi_i)\sum_{k=0}^{\infty}\frac{\Pi_{i,k+1}^{a}(T^{n}\xi_i)}{\Pi_{i,k}^{ p-1 }\left(T^{n}\xi_i\right)}\left(\frac{m_{i,k}^{(p)}\left(T^{n}\xi_i\right)}{m_{i,k}^{p}\left(T^{n}\xi_i\right)}\right) \right].
\end{align*}
Since $\Pi_{i,n}(\xi_i)  $ is independent of $ \Pi_{i,k+1}(T^{n}\xi_i) $ and $ \Pi_{i,k}(T^{n}\xi_i) $ under $ \mathbb{P} $, for any $ k\geq 0 $,
\begin{align}
	q_{i,n,K}\mathbb{E}\widetilde{N}_{i,n,K}^{-a}
	&=\mathbb{E}\left[\Pi_{i,n}^{a}\,\prod_{j=0}^{n-1}\Big(p_{1}(\xi_{i,j})+(1-p_{1}(\xi_{i,j}))\beta_{K}\Big)\right]\nonumber
	\\
	& \quad + \frac{1}{K }\,\mathbb{E}\Pi_{i,n}^{a}\,\mathbb{E}\left[1+\sum_{k=0}^{\infty}\frac{\Pi_{i,k+1}^{a}(T^{n}\xi_i)}{\Pi_{i,k}^{ p-1  }\left(T^{n}\xi_i\right)}\left(\frac{m_{i,k}^{(p)}\left(T^{n}\xi_i\right)}{m_{i,k}^{p}\left(T^{n}\xi_i\right)}\right) \right]\nonumber\\
	&=\mathbb{E}\prod_{j=0}^{n-1}\left[m_{i,j}^{a}\left(p_{1}(\xi_{i,j})+(1-p_{1}(\xi_{i,j}))\beta_{K}\right)\right]\nonumber
	\\
	&\quad + \frac{1}{K }\mathbb{E}\Pi_{i,n}^{a}\,\mathbb{E}\left[1+\sum_{k=0}^{\infty}\frac{\Pi_{i,k+1}^{a}(T^{n}\xi_i)}{\Pi_{i,k}^{ p-1  }\left(T^{n}\xi_i\right)}\left(\frac{m_{i,k}^{(p)}\left(T^{n}\xi_i\right)}{m_{i,k}^{p}\left(T^{n}\xi_i\right)}\right) \right]\nonumber\\
	&= \Big\{\mathbb{E}\left[m_{i,0}^{a}\left(p_{1}(\xi_{i,0})+(1-p_{1}(\xi_{i,0}))\beta_{K}\right)\right] \Big \}^{n} \label{v5}\\
	&\quad + \frac{(\mathbb{E}m_{i,0}^{a} )^n }{K } \left[1 +\frac{1}{1 - \mathbb{E} m_{i,0}^{a-(p-1)}}  \mathbb{E}\Big(\frac{ m_{i,0}^{a}m_{i,0}^{(p)} }{m_{i,0}^{p} }\Big)\right].\nonumber
\end{align}
Notice that $0\leq  m_{i,0}^{a}\left(p_{1}(\xi_{i,0})+(1-p_{1}(\xi_{i,0}))\beta_{K}\right)  \leq m_{i,0}^{a}$. By condition \textbf{A3}, it holds $\mathbb{E} m_{i,0}^{\lambda_0} < \infty.$ Thus by the dominated convergence theorem, we have
$$\lim_{a\searrow 0} \Big\{\mathbb{E}\left[m_{i,0}^{a}\left(p_{1}(\xi_{i,0})+(1-p_{1}(\xi_{i,0}))\beta_{K}\right)\right] \Big \}^{n} = \Big\{\mathbb{E}\left[  p_{1}(\xi_{i,0})+(1-p_{1}(\xi_{i,0}))\beta_{K} \right] \Big \}^{n},$$
and
$$ \lim_{a\searrow 0}(\mathbb{E}m_{i,0}^{a} )^n =1 . $$
Since $m_{i,0} \geq 1$ a.s. and $p >1$, we have
$ \lim_{a\searrow 0}\mathbb{E} m_{i,0}^{a-(p-1)}=\mathbb{E} m_{i,0}^{-(p-1)}.$
When $a \in (0, p-1],$ we have  $ 0\leq \frac{ m_{i,0}^{a}m_{i,0}^{(p)} }{m_{i,0}^{p} }  \leq   \frac{  m_{i,0}^{(p)} }{m_{i,0}  }.  $
Condition \textbf{A4} implies that $$ \mathbb{E}  \frac{  m_{i,0}^{(p)} }{m_{i,0}  }=   \mathbb{E}  \frac{ \mathbb{E} [ Z_{i,1}^{p}| \xi_0] }{m_{i,0}  }  =  \mathbb{E}  \frac{Z_{1,1}^{p}}{m_{1,0}}  < \infty. $$
 Thus by the dominated convergence theorem, we have $$ \lim_{a\searrow 0} \mathbb{E}\Big(\frac{ m_{i,0}^{a}m_{i,0}^{(p)} }{m_{i,0}^{p} }\Big)= \mathbb{E}\Big(\frac{  m_{i,0}^{(p)} }{m_{i,0}^{p} }\Big).$$ Hence, from (\ref{v5}), we get
 \begin{eqnarray*}
  && q_{i,n,K}\mathbb{E}\widetilde{N}_{i,n,K}^{-a}\stackrel{a\downarrow0}{\longrightarrow}q_{i,n,K} = \Big\{\mathbb{E}\left[  p_{1}(\xi_{i,0})+(1-p_{1}(\xi_{i,0}))\beta_{K} \right] \Big \}^{n}  \\
	&&\ \ \ \ \ \ \ \ \  \ \ \ \ \ \ \ \ \ \ \ \ \ \ \ \ \ \ \ \ \ \ \ \ \ \ \  +\ \frac{1 }{K } \left[1 +\frac{1}{1 - \mathbb{E} m_{i,0}^{-(p-1)}}  \mathbb{E}\Big(\frac{  m_{i,0}^{(p)} }{m_{i,0}^{p} }\Big)\right].
\end{eqnarray*}
As $\beta_{K} \in (0, 1),$ we have $\mathbb{E}\left[  p_{1}(\xi_{i,0})+(1-p_{1}(\xi_{i,0}))\beta_{K} \right] < 1$, which leads to
$$\lim_{n \rightarrow \infty}\Big\{\mathbb{E}\left[  p_{1}(\xi_{i,0})+(1-p_{1}(\xi_{i,0}))\beta_{K} \right] \Big \}^{n}=0.$$
Therefore, it holds
 \begin{eqnarray*}
  &&  q_{i,n,K}\stackrel{n \rightarrow \infty}{\longrightarrow}   \frac{1 }{K } \left[1 +\frac{1}{1 - \mathbb{E} m_{i,0}^{-(p-1)}}  \mathbb{E}\Big(\frac{  m_{i,0}^{(p)} }{m_{i,0}^{p} }\Big)\right]\stackrel{K\rightarrow\infty}{\longrightarrow}0.
\end{eqnarray*}
In conclusion, we have
\begin{equation*} q_{i,n,K}\mathbb{E}\widetilde{N}_{i,n,K}^{-a}\stackrel{a\downarrow0}{\longrightarrow}q_{i,n,K}\stackrel{n\rightarrow\infty}{\longrightarrow}\frac{1}{K }\left(1+C\right)\stackrel{K\rightarrow\infty}{\longrightarrow}0.
\end{equation*}
Then, we take $n_0, K_{0}$ and $ a_0 \in (0, \lambda_0)$ small enough such that $q_{i,n_0,K_0}\mathbb{E} \widetilde{N}_{i,n_0,K_0}^{-a_0} < 1$.
By \eqref{vs4} and Markov's inequality, we can get
\begin{align}
\phi_i(t)&\leq q_{i,n_0,K_0}\mathbb{E}\phi(\widetilde{N}_{i,n_0,K_0}\,t)+\mathbb{P}\Big(\frac{t}{\Pi_{i,n_0}}<t_{K_{0}}\Big)\nonumber \\
&\leq q_{i,n_0,K_0}\mathbb{E}\phi(\widetilde{N}_{i,n_0,K_0}\,t)+ \frac{C_{n_0, K_0}}{t^{a_0}} .\label{l65}
\end{align}
Notice that $0\leq  \phi_i(t)\leq 1\ (t>0)$.
Finally, by \eqref{l65} and Lemma 4.1 in \cite{Liu1999}, we have for all $t > 0$,
\begin{equation*}
	\phi_i(t)\leq \frac{C}{ 1+ t^{a_{0}}} .
\end{equation*}
Let $0<\alpha<a_{0}$, we have
\begin{equation}\label{i4.47}
	\int_{1}^{\infty} \phi_i(t) t^{\alpha-1} dt\leq C \int_{1}^{\infty}  t^{\alpha-a_{0}-1} dt<\infty.
\end{equation}
By \eqref{i4.45}, \eqref{i4.46} and \eqref{i4.47},  inequality \eqref{eqLidsu3} holds.

It remains to prove \eqref{eqLisu4} now.   Since $x\mapsto x^{-\alpha}\ (\alpha>0,\ x>0)$ is a non-negative convex function.
Then by Lemma 2.1 in \cite{HL12}, we have
$$\sup_{n\in\mathbb{N}}\mathbb{E}  W_{i, n}^{-\alpha}=\mathbb{E}  W_{i, \infty}^{-\alpha}<\infty.$$
This completes the proof of Lemma \ref{lemma1}.
\hfill\qed

Under the conditions \textbf{A3} and \textbf{A4}, we have the following analogue to  Lemma \ref{lemma4.5}.

 \begin{lemma}\label{lama4.9}
 		Assume that the conditions \textbf{A3} and \textbf{A4} are satisfied. Then  for  all  $   |x|  \leq  \sqrt{\ln (m \wedge n)} ,$
 		\begin{equation}  \label{ff7}
 			\mathbb{P}\bigg(R_{m,n} \leq
 			x, \sum_{i=1}^{n+m}  \eta_{ m,n, i}\geq x\bigg) \leq  C  \frac{1+ x^2 }{ \sqrt{m \wedge n} \   }  \exp\bigg\{ - \frac12 x^2   \bigg\}
 		\end{equation}
 		and
 		\begin{equation} \label{ff11}
 			\mathbb{P}\bigg(R_{m,n} \geq x, \sum_{i=1}^{n+m}  \eta_{ m,n, i}\leq x\bigg) \leq
 			C  \frac{1+ x^2 }{ \sqrt{m \wedge n} \   } \exp\bigg\{ - \frac12 x^2   \bigg\}.
 		\end{equation}		
 \end{lemma}
\emph{Proof.}  As the conditions \textbf{A3} and \textbf{A4} together imply the conditions \textbf{A1} and \textbf{A2}, when  $  |x|  \leq 1$,  the inequalities \eqref{ff7} and \eqref{ff11} are simple consequences of Lemma \ref{lemma4.5}.
Thus we only need to prove the inequalities \eqref{ff7} and \eqref{ff11}  for  all  $ 1 \leq |x|  \leq  \sqrt{\ln (m \wedge n)} .$
In the sequel, we only give a proof for inequality (\ref{ff7}) with $ 1 \leq |x|  \leq  \sqrt{\ln (m \wedge n)}$, as inequality  \eqref{ff11} with  $ 1 \leq |x|  \leq  \sqrt{\ln (m \wedge n)} $ can be proved in the same way.

 Without loss of generality, we may assume that $  m\leq n$.  Recall  the  notations $Y_{m,n,k},$ $\widetilde{Y}_{m,n, k}$, $H_{m,n,k}$ and $D_{m, n , k}$ defined by (\ref{Y1}), (\ref{Y2}), (\ref{Y3}) and (\ref{Y4}), respectively. Set $ \alpha_{m}= m^{-1/2}$ and $k=[ m^{1/2} \,]$.
From \eqref{4.38}, it holds for all $  x \in \mathbb{R}$,
\begin{align}
&\mathbb{P}\Big(R_{m,n} \leq x, \sum_{i=1}^{n+m}  \eta_{ m,n, i}\geq x\Big)
\nonumber\\
\leq\ & \mathbb{P}\Big(Y_{m,n,0}+H_{m,n,k}\leq x+ \alpha_{m}, Y_{m,n,0}\geq x\Big)  +\ \mathbb{P}\Big(|D_{m, n , k}|\geq \alpha_{m}\Big). \label{4.d38}
\end{align}
For $\mathbb{P}\left(|D_{m, n , k}|\geq \alpha_{m}\right)$,   by an argument similar to (\ref{4.39}), we have  for all  $ 1 \leq |x| \leq \sqrt{\ln m}$,
\begin{equation}\label{4.d39}
	\mathbb{P}\left(\left|D_{m, n, k}\right|>\alpha_{m}\right)\leq C \frac{x^2 }{ \sqrt{m  } \   }  \exp\Big\{ - \frac12 x^2   \Big\}.
\end{equation}
Next, we give an estimation for the first term of  bound \eqref{4.d38}. Recall the notations $G_{m, n,k}(x)$, $v_{k}(d s, d t)$
defined by (\ref{inds25a}). Recall $C_{m,n,k}^2 =  \textrm{Var} (Y_{m,n,k})$. Then we have
\begin{align*}
	&\mathbb{P}\left(Y_{m,n,0}+H_{m,n,k} \leq x+\alpha_{m},  Y_{m,n,0} \geq x\right)\nonumber \\
 =\ &  \iint \mathbf{1}_{\{t \leq \alpha_{m}\}}\Big(G_{m,  n,k}\left(x-s-t+\alpha_{m}\right)-G_{m,  n,k}(x-s)\Big) v_{k}(ds,  dt)
\end{align*}
and $ C_{m,n,k}=  1 +O(1/\sqrt{m}) \rightarrow 1  $ as $ m  \rightarrow \infty.$
Using Cram\'{e}r's moderate deviations (for $|x| \leq m^{1/6}$) and Bernstein's inequality (for $|x| > m^{1/6}$) for independent random variables, we have the following  non-uniform Berry-Esseen's bound: for all $   x \in \mathbb{R} ,$
\begin{align*}
&\left|\mathbb{P}\left(  \frac{ Y_{m,n,k} }{C_{m,n,k} }  \leq \frac{x}{C_{m,n,k}}\right)-\Phi\left(\frac{x}{C_{m,n,k}}\right)\right|\\
\leq\ & C_1 \exp\bigg\{  - \frac{x^2 }{ 2(1+ \frac{C}{\sqrt{m}} |x|) }   \bigg\}\Bigg( 1+ \Big(\frac{|x|}{C_{m,n,k} }  \Big)^2   \Bigg) \\
&\quad  \times \Bigg(\sum_{j=k+1}^{m}\mathbb{E}\left|\frac{\eta_{ m,n, i}}{ C_{m,n,k}} + \frac{\eta_{ m,n, n+ i}}{ C_{m,n,k}}\right|^{3} + \sum_{i=m+1}^{n}\mathbb{E}\left|\frac{\eta_{ m,n, i}}{ C_{m,n,k}}\right|^{3}
\Bigg)  \\
\leq\ &    C_2 \frac{ 1+ |x|^2 }{\sqrt{m} }  \exp\bigg\{  - \frac{x^2 }{ 2 \big( 1+ \frac{C}{\sqrt{m}} |x| \big) }   \bigg\}.
\end{align*}
By the last inequality, we deduce that for all  $   x \in \mathbb{R} ,$
\begin{align*}
&\left|G_{m,  n,k}(x)-\Phi(x)\right| \\
\leq\ & \left|\mathbb{P}\left( \frac{ Y_{m,n,k} }{C_{m,n,k} }\leq \frac{x}{C_{m,n,k}}\right)-\Phi\left(\frac{x}{C_{m,n,k}}\right)\right|
  +  \left|\Phi\left(\frac{x}{C_{m,n,k}}\right)-\Phi\left(x\right)\right|\\
\leq\ & C_2 \frac{ 1+ x^2 }{\sqrt{m} }  \exp\bigg\{  - \frac{x^2 }{  2 \big( 1+ \frac{C}{\sqrt{m}} |x| \big) }   \bigg\} +\exp\bigg\{-\frac{x^{2}}{2 }\bigg\}\left|\frac{x}{C_{m,n,k}}-x\right|\\
\leq\ &   C_2 \frac{ 1+ x^2 }{\sqrt{m} }  \exp\bigg\{  - \frac{x^2 }{  2 \big( 1+ \frac{C}{\sqrt{m}} |x| \big) }   \bigg\} +  C_3 \frac{|x|}{\sqrt{m}}   \exp\bigg\{-\frac{x^{2}}{2  }\bigg\} \\
\leq\ &  C_4 \frac{ 1+ x^2 }{\sqrt{m} }  \exp\bigg\{  - \frac{x^2 }{  2 \big( 1+ \frac{C}{\sqrt{m}} |x| \big) }   \bigg\} .
\end{align*}
Therefore,   we have for all $x \in \mathbb{R},$
\begin{equation}\label{4.d43}
  \mathbb{P}\Big(Y_{m,n,0}+H_{m,n,k} \leq x+\alpha_{m},  Y_{m,n,0} \geq  x\Big) \leq  J_{1}+J_{2} +J_3,
\end{equation}
where
$$J_{1}=\iint \mathbf{1}_{\{t \leq \alpha_{m}\}}\Big|\Phi\left(x-s-t+\alpha_{m}\right)-\Phi(x-s)\Big| v_{k}(d s,  d t),$$
$$ J_{2}= C \iint \mathbf{1}_{\{t \leq \alpha_{m}\}} \frac{ 1+ |x-s|^2 }{\sqrt{m} }  \exp\bigg\{  - \frac{(x-s)^2 }{ 2 \big( 1+ \frac{C}{\sqrt{m}} |x-s| \big) }   \bigg\}   v_{k}(d s,  d t)$$
and
$$ J_{3}= C \iint \mathbf{1}_{\{t \leq \alpha_{m}\}} \frac{ 1+ |x-s-t|^2 }{\sqrt{m} }  \exp\bigg\{  - \frac{(x-s-t)^2 }{
 2 \big(1+ \frac{C}{\sqrt{m}} |x-s-t| \big)}   \bigg\}  v_{k}(d s,  d t) . $$
Denote $\widetilde{C}_{m,n,k}^2= \textrm{Var}(\widetilde{Y}_{m,n, k} )$, then it holds $\widetilde{C}_{m,n,k}^2=O(1/\sqrt{m})$ as $m\rightarrow \infty$.
For the upper bound of $J_{1}$, by the mean value theorem, we have for all $1\leq |x|   \leq \sqrt{\ln m},$
\begin{align*}
 &\mathbf{1}_{\{t \leq \alpha_{m}\}}\left|\Phi\left(x-s-t+\alpha_{m}\right)-\Phi(x-s)\right| \\
 \leq \ & C|\alpha_{m}-t |\exp\bigg\{  - \frac{x^2 }{  2 \big( 1+ \frac{C}{\sqrt{m}} |x| \big) }   \bigg\} + \Big|\alpha_{m}-t  \Big| \mathbf{1}_{\{|s| \geq \,  2|x| \widetilde{C}_{m,n,k}    \}} +   \,  \Big|\alpha_{m}-t \Big| \mathbf{1}_{\{|t| \geq  \, C_0  |x|  \widetilde{C}_{m,n,k}   \}},
\end{align*}
which leads to
\begin{equation}\label{4.d44}
	J_{1}\leq J_{11}+ J_{12}+ J_{13},
\end{equation}
where
\begin{align*}
	J_{11}&= C\iint |\alpha_{m}-t |\exp\bigg\{  - \frac{x^2 }{  2 \big( 1+ \frac{C}{\sqrt{m}} |x| \big) }   \bigg\} v_{k}(d s,  d t),  \\
J_{12}&=  \iint |\alpha_{m}-t | \mathbf{1}_{\{|s| \geq  \, 2 |x| \widetilde{C}_{m,n,k}   \}}  v_{k}(d s,  d t)
\end{align*}
and
\begin{equation*}
	J_{13}=  \iint |\alpha_{m}-t | \mathbf{1}_{\{|t| \geq \, C_0  |x|  \widetilde{C}_{m,n,k} \}} v_{k}(d s,  d t).
\end{equation*}
By Lemma \ref{lemma4.3}, it is obvious that  for all   $1\leq |x|   \leq \sqrt{\ln m},$
\begin{align*}
	J_{11}& \leq   C_1 \bigg(\alpha_{m} + \mathbb{E}| H_{m,n,k}| \bigg)\exp\bigg\{  - \frac{x^2 }{  2 \big( 1+ \frac{C}{\sqrt{m}} |x| \big) }   \bigg\} \\
&\leq    \frac{C_2}{\sqrt{m} } \exp\bigg\{  - \frac{x^2 }{  2 \big( 1+ \frac{C}{\sqrt{m}} |x| \big) }   \bigg\}.
\end{align*}
For $J_{12},$ we have the following estimation for all   $1\leq |x|   \leq \sqrt{\ln m},$
\begin{equation*}
	J_{12}\leq       \alpha_{m} \mathbb{P}\bigg( |\widetilde{Y}_{m,n, k} | \geq  2 |x| \widetilde{C}_{m,n,k}   \bigg) + \mathbb{E}| H_{m,n,k}|\textbf{1}_{\{|\widetilde{Y}_{m,n, k} | \geq 2 |x| \widetilde{C}_{m,n,k} \}}   .
\end{equation*}
By Bernstein's inequality, we deduce that for all $x \in \mathbb{R},$
\begin{align*}
 \mathbb{P}\bigg( |\widetilde{Y}_{m,n, k} | \geq   2 |x| \widetilde{C}_{m,n,k}    \bigg)
&=\mathbb{P}\left(  \frac{\widetilde{Y}_{m,n, k} }{\widetilde{C}_{m,n,k} } \geq 2  |x| \right)+\mathbb{P}\left(   \frac{\widetilde{Y}_{m,n, k} }{\widetilde{C}_{m,n,k} }  \leq - 2  |x|\right)\\
&\leq 2 \exp\bigg\{  - \frac{(2x)^2 }{  2 \big( 1+ \frac{C}{\sqrt{m}} |x| \big) }   \bigg\} . \nonumber
\end{align*}
By Cauchy-Schwartz's inequality, we have for all $x \in \mathbb{R},$
\begin{align*}
\mathbb{E}| H_{m,n,k}|\textbf{1}_{\{|\widetilde{Y}_{m,n, k} | \geq  2 |x| \widetilde{C}_{m,n,k} \}}  & \leq
\Big(\mathbb{E}| H_{m,n,k}|^2 \Big)^{1/2}  \mathbb{P}\bigg( |\widetilde{Y}_{m,n, k} | \geq   2 |x| \widetilde{C}_{m,n,k}    \bigg) ^{1/2}     \\
  & \leq  \frac{C}{\sqrt{m} }   \Bigg(2\exp\bigg\{  - \frac{(2x)^2 }{  2 \big( 1+ \frac{C}{\sqrt{m}} |x| \big) }   \bigg\} \Bigg)^{1/2}  \nonumber  \\
  & \leq   \frac{C_1}{\sqrt{m} }  \exp\Bigg\{  - \frac{x^2 }{  2 \big( 1+ \frac{C}{\sqrt{m}} |x| \big) }   \Bigg\}.
\end{align*}
Hence, we get for  all $1\leq|x|  \leq \sqrt{ \ln m } ,$
\begin{equation*}
	J_{12} \leq   \frac{C}{\sqrt{m}\ \ } \exp\bigg\{  - \frac{x^2 }{  2 \big( 1+ \frac{C}{\sqrt{m}} |x| \big) }   \bigg\}  .
\end{equation*}
For $J_{13},$  we have   for all   $1\leq |x|   \leq \sqrt{\ln m},$
\begin{equation*}
	J_{13} \leq     \alpha_{m} \mathbb{P}\bigg( |H_{m,n,k}| \geq C_0  |x| \widetilde{C}_{m,n,k}  \bigg) +    \mathbb{E} \Big[| H_{m,n,k}| \textbf{1}_{\{|H_{m,n,k}| \geq C_0  |x|  \widetilde{C}_{m,n,k} \} } \Big]    .  \nonumber
\end{equation*}
Notice that $V_{m,n,\rho}\asymp \frac{ 1}{\sqrt{m} }$ and $\widetilde{C}_{m,n,k}  \asymp \frac{ 1}{ m^{1/4} }$.
It is easy to see  that for all   $1\leq |x|   \leq \sqrt{\ln m},$
\begin{equation*}
 \mathbb{P}\bigg( |H_{m,n,k} | \geq C_0 |x|  \widetilde{C}_{m,n,k}  \bigg)
 \leq T_1 + T_2,
\end{equation*}
where
\begin{equation*}
 T_1 =  \mathbb{P}\bigg( \Big|\frac{\ln W_{1, k}}{ n\, V_{m,n,\rho}  \ } \Big|    \geq \frac12 C_0 |x| \widetilde{C}_{m,n,k}  \bigg) \
    \ \textrm{and} \ \     \ T_2=\mathbb{P}\bigg(   \Big| \frac{\ln W_{2,k}}{ m\, V_{m,n,\rho} \ } \Big| \geq \frac12 C_0 |x|  \widetilde{C}_{m,n,k}  \bigg)  .
\end{equation*}
By Lemma \ref{lemma1} and  Markov's inequality,  we have for all   $1\leq |x|   \leq \sqrt{\ln m},$
\begin{align*}
T_1 &\leq  \mathbb{P}\bigg( \frac{\ln W_{1, k}}{ n\, V_{m,n,\rho}  \ }   \geq \frac12 C_0 |x| \widetilde{C}_{m,n,k}  \bigg) + \mathbb{P}\bigg( \frac{\ln W_{1, k}}{ n\, V_{m,n,\rho}  \ }  \leq - \frac12 C_0 |x| \widetilde{C}_{m,n,k}  \bigg)  \\
&\leq   \mathbb{P}\bigg( W_{1, k}    \geq \exp\{ \frac12 C_0 |x| n\, V_{m,n,\rho} \widetilde{C}_{m,n,k}\}  \bigg) + \  \mathbb{P}\bigg( W_{1, k}^{-1}    \geq \exp\{ \frac12 C_0 |x| n\, V_{m,n,\rho} \widetilde{C}_{m,n,k}\}  \bigg)  \\
&\leq  \mathbb{E}[W_{1,k} ]  \exp\bigg\{ -\frac12 C_0 |x| n\, V_{m,n,\rho} \widetilde{C}_{m,n,k} \bigg\}  + \ \mathbb{E}[W_{1,k}^{-\alpha} ]  \exp\bigg\{ -\frac12 \alpha C_0 |x| n\, V_{m,n,\rho} \widetilde{C}_{m,n,k} \bigg\}  \\
& \leq   C   \exp\bigg\{ - \frac12 x^2 \bigg\} ,
\end{align*}
with $C_0$ large  enough. Similarly, we have for all   $1\leq |x|   \leq \sqrt{\ln m},$
\begin{equation*}
T_2 \leq    C   \exp\bigg\{ - \frac12 x^2 \bigg\} . \nonumber
\end{equation*}
Hence, we get   for all   $1\leq |x|   \leq \sqrt{\ln m},$
\begin{equation*}
 \mathbb{P}\bigg( |H_{m,n,k} | \geq C_0 |x|  \widetilde{C}_{m,n,k}  \bigg)   \leq    C   \exp\bigg\{ - \frac12 x^2 \bigg\} .
\end{equation*}
Clearly, by Lemma \ref{lemma1} and the inequality $|\ln x|^2 \leq C_\alpha (x+ x^{-\alpha})$ for all $\alpha, x>0$, it holds
 $$\mathbb{E}| H_{m,n,k}|^2  \leq \frac{C_1}{m}\Big(\mathbb{E} W_{1, k} +\mathbb{E} W_{1, k}^{-\alpha} + \mathbb{E} W_{2, k}+\mathbb{E} W_{2, k}^{-\alpha}  \Big) \leq  \frac{C_2}{m}.$$
By Markov's inequality, Cauchy-Schwartz's inequality and Lemma \ref{lemma1}, we deduce that   for all   $1\leq |x|   \leq \sqrt{\ln m},$
\begin{align*}
& \mathbb{E}\Big[| H_{m,n,k}| \textbf{1}_{\{|H_{m,n,k}|   \geq C_0 | x|  \widetilde{C}_{m,n,k} \} }\Big]\\
  \leq\ & \mathbb{E}\Big[| H_{m,n,k}| \textbf{1}_{\{W_{1, k}    \geq \exp\{ \frac12 C_0 |x|  n\, V_{m,n,\rho} \widetilde{C}_{m,n,k}\}  \} }  \Big]\\
  &\quad +\mathbb{E}\Big[| H_{m,n,k}| \textbf{1}_{\{W_{2,k} \geq \exp\{ \frac12  C_0 |x|  m\, V_{m,n,\rho}\widetilde{C}_{m,n,k} \} \} } \Big] \\
&\quad+  \mathbb{E}\Big[| H_{m,n,k}| \textbf{1}_{\{W_{1, k}^{-1}    \geq \exp\{ \frac12 C_0 |x|  n\, V_{m,n,\rho} \widetilde{C}_{m,n,k}\}  \} }  \Big]\\
&\quad +   \mathbb{E}\Big[| H_{m,n,k}| \textbf{1}_{\{W_{2,k}^{-1}  \geq \exp\{ \frac12  C_0 |x|  m\, V_{m,n,\rho}\widetilde{C}_{m,n,k} \} \} } \Big]  \\
  \leq& \exp\bigg\{ -  \frac 1  4 C_0 |x| n\, V_{m,n,\rho} \widetilde{C}_{m,n,k} \bigg\}    \mathbb{E}[W_{1, k}^{1/2} | H_{m,n,k}|]\\
  &\quad +\exp\bigg\{ -  \frac 1  4 C_0 |x| m\, V_{m,n,\rho} \widetilde{C}_{m,n,k} \bigg\}    \mathbb{E}[W_{2, k}^{1/2} | H_{m,n,k}|] \\
 &\quad + \ \exp\bigg\{  - \frac 1  4 \alpha C_0 |x|  n\, V_{m,n,\rho}\widetilde{C}_{m,n,k} \bigg\}\mathbb{E}[W_{1,k}^{-\alpha/2}  | H_{m,n,k}|  ]\\
 &\quad + \ \exp\bigg\{  - \frac 1  4 \alpha C_0 |x|  m\, V_{m,n,\rho}\widetilde{C}_{m,n,k} \bigg\}\mathbb{E}[W_{2,k}^{-\alpha/2}  | H_{m,n,k}|  ]\\
    \leq\  &\exp\bigg\{ -  \frac 1  4 C_0 |x| n\, V_{m,n,\rho} \widetilde{C}_{m,n,k} \bigg\}   ( \mathbb{E}W_{1, k})^{1/2} ( \mathbb{E}| H_{m,n,k}|^2 )^{1/2}\\
    &\quad  + \exp\bigg\{ -  \frac 1  4 C_0 |x| m \, V_{m,n,\rho} \widetilde{C}_{m,n,k} \bigg\}   ( \mathbb{E}W_{2, k})^{1/2} ( \mathbb{E}| H_{m,n,k}|^2 )^{1/2}   \\
    &\quad  + \exp\bigg\{ -  \frac 1  4 \alpha C_0 |x| n\, V_{m,n,\rho} \widetilde{C}_{m,n,k} \bigg\}   ( \mathbb{E}W_{1, k}^{-\alpha })^{1/2} ( \mathbb{E}| H_{m,n,k}|^2 )^{1/2}\\
    &\quad + \exp\bigg\{  - \frac 1  4 \alpha C_0 |x|  m\, V_{m,n,\rho}\widetilde{C}_{m,n,k} \bigg\}  \ (\mathbb{E} W_{2,k}^{-\alpha } )^{1/2}( \mathbb{E}| H_{m,n,k}|^2 )^{1/2} \\
   \leq\ & \frac{C}{\sqrt{m} } \exp\bigg\{ - \frac12 x^2 \bigg\} ,
\end{align*}
where the last inequality holds with $C_0$ large  enough.
Hence, we have    for all   $1\leq |x|   \leq \sqrt{\ln m},$
\begin{equation*}
	J_{13}\leq   \frac{C}{\sqrt{m} } \exp\Big\{ - \frac12 x^2 \Big\} .
\end{equation*}
Returning to (\ref{4.d44}), we get   for all   $1\leq |x|   \leq \sqrt{\ln m},$
\begin{equation} \label{fgssgd43}
	J_{1} \leq  \frac{C}{\sqrt{m} }\exp\bigg\{  - \frac{x^2 }{  2 \big( 1+ \frac{C}{\sqrt{m}} |x| \big) }   \bigg\} .
\end{equation}
For $J_{2}$, by an argument similar to the proof of \eqref{fgssgd43}, we have for all   $1\leq |x|   \leq \sqrt{\ln m},$
\begin{align}
	J_{2}&=C_1 \iint \mathbf{1}_{\{t \leq \alpha_{m}\}} \frac{ 1+ |x-s|^2 }{\sqrt{m} }  \exp\bigg\{  - \frac{(x-s)^2 }{ 2 \big( 1+ \frac{C_2}{\sqrt{m}} |x-s| \big) }   \bigg\}   v_{k}(d s,  d t)\nonumber  \\
  &\leq\frac{C_3}{\sqrt{m} }\left(\int_{|s|\leq |x|\widetilde{C}_{m,n,k} } (1+x^2 ) \exp\bigg\{  - \frac{x^2 }{ 2 \big( 1+ \frac{C_4}{\sqrt{m}} |x | \big) }   \bigg\} v_{k}(d s) +\int_{|s| > |x| \widetilde{C}_{m,n,k} } (1 + x^2 )v_{k}(d s)\right)  \nonumber \\
	&\leq    \frac{C_3}{\sqrt{m} }(1+x^2) \left[ \exp\bigg\{  - \frac{x^2 }{ 2 \big( 1+ \frac{C_4}{\sqrt{m}} |x | \big) }   \bigg\}+\mathbb{P}\bigg( \Big| \frac{\widetilde{Y}_{m,n,k} }{\widetilde{C}_{m,n,k} } \Big|  > |x| \bigg)\right]\nonumber  \\
	&\leq    \frac{C }{\sqrt{m} } (1+x^2) \exp\bigg\{  - \frac{x^2 }{  2 \big( 1+ \frac{C}{\sqrt{m}} |x| \big) }   \bigg\}.\label{fgsgd44}
\end{align}
Similarly, for $J_{3}$,  we have  for all $1\leq |x|\leq \sqrt{\ln m},$
\begin{align}
J_{3}&= C_1 \iint  \mathbf{1}_{\{t \leq \alpha_{m}\}} \frac{ 1+ |x-s-t|^2 }{\sqrt{m} }  \exp\bigg\{  - \frac{(x-s-t)^2 }{
 2 \big(1+ \frac{C}{\sqrt{m}} |x-s-t| \big)}   \bigg\}  v_{k}(d s,  d t) \nonumber  \\
&\leq \frac{C_2}{\sqrt{m} }\left(\iint  (1+x^2) \exp\bigg\{  - \frac{x^2 }{
 2 \big(1+ \frac{C}{\sqrt{m}} |x| \big)}   \bigg\}  v_{k}(d s,  d t)  \right. \nonumber \\
 &\quad +\iint_{|s| > |x|\widetilde{C}_{m,n,k}} v_{k}(d s,  d t) + \left.  \iint_{|t| >C_0 |x| \widetilde{C}_{m,n,k}} v_{k}(d s,  d t) \right)  \nonumber \\
	&\leq \frac{C_2}{\sqrt{m} }\left(  (1+x^2) \exp\bigg\{  - \frac{x^2 }{
 2 \big(1+ \frac{C}{\sqrt{m}} |x| \big)}   \bigg\}+\mathbb{P}\bigg(  \Big| \frac{\widetilde{Y}_{m,n,k} }{\widetilde{C}_{m,n,k} } \Big|  > |x|     \bigg)   \right.  \nonumber  \\
 &\quad +  \mathbb{P}\bigg(   |H_{m,n,k}|  > C_0 |x| \widetilde{C}_{m,n,k}   \bigg) \Bigg)\nonumber  \\
	 &\leq \frac{C_4}{\sqrt{m} } (1+x^2) \exp\bigg\{  - \frac{x^2 }{  2 \big( 1+ \frac{C}{\sqrt{m}} |x| \big) }   \bigg\}. \label{fgsgd45}
\end{align}
Applying  the inequalities \eqref{fgssgd43}-\eqref{fgsgd45}  to \eqref{4.d43}, we get  for all   $1\leq |x|   \leq \sqrt{\ln m},$
\begin{align}
&\mathbb{P}\Big(Y_{m,n,0}+H_{m,n,k}\leq x+ \alpha_{m}, Y_{m,n,0}\geq x\Big) \nonumber \\
\leq\ &   \frac{C_1}{\sqrt{m} } (1+x^2) \exp\bigg\{  - \frac{x^2 }{  2 \big( 1+ \frac{C}{\sqrt{m}} |x| \big) }   \bigg\} \nonumber \\
\leq\ &   \frac{C_1}{\sqrt{m} } (1+x^2)(1 +  \frac{C_2}{\sqrt{m}} |x|^3 ) \exp\bigg\{  - \frac{x^2 }{  2 }     \bigg\} \nonumber \\
\leq\ &   \frac{C_3}{\sqrt{m} } (1+x^2) \exp\bigg\{  - \frac{x^2 }{  2 }     \bigg\}.\label{u4.d40}
\end{align}
Combining  (\ref{4.d38}),  (\ref{4.d39}) and \eqref{u4.d40} together, we get  (\ref{ff7})   for all   $1\leq |x|   \leq \sqrt{\ln m}.$
This completes the proof of Lemma \ref{lama4.9}. \hfill\qed

\subsection{Proof of Theorem \ref{th2.2} }   \label{sec5}

We give a  proof of Theorem \ref{th2.2} for the case of $\frac{\mathbb{P} ( R_{m,n} \geq x  )}{1-\Phi(x)},$ $ x\geq 0.$ Thanks to the
symmetry between $m$ and $n$, the case  of $\frac{\mathbb{P} (  -R_{m,n}  \geq x  )}{ \Phi(-x)}$ can be proved in the similar way.
To prove Theorem \ref{th2.2}, we start with the proofs of Lemmas   \ref{th200}  and  \ref{th100}, and conclude
Theorem \ref{th2.2} by combining  Lemmas   \ref{th200}  and  \ref{th100} together.
To avoid  trivial  case,  we assume that $ m \wedge n \geq 2$.

The following lemma gives the upper bound in Theorem \ref{th2.2}.
\begin{lemma} \label{th200}
Assume that the conditions \textbf{A3} and \textbf{A4} are satisfied. Then  it holds
 for   all  $0 \leq  x \leq  c \, \sqrt{m \wedge n}  ,$
\begin{equation}\label{ineq332}
\ln  \frac{\mathbb{P}\big(  R_{m,n}  \geq x  \big)}{1-\Phi(x)}  \leq C  \frac{ 1+ x^3   }{  \sqrt{m \wedge n}\  } .
\end{equation}
\end{lemma}
\emph{Proof.}  First, we consider the   case  $0 \leq x  \leq  \sqrt{\ln(m \wedge n)} $.
Notice that
\begin{align*}
	\mathbb{P}\Big(R_{m,n}\geq x\Big)
	&= \mathbb{P}\Big(R_{m,n} \geq x, \sum_{i=1}^{n+m}  \eta_{ m,n, i}  \geq x\Big) + \mathbb{P}\Big( R_{m,n} \geq x, \sum_{i=1}^{n+m}  \eta_{ m,n, i}  < x\Big) \\
	&\leq \mathbb{P}\Big( \sum_{i=1}^{n+m}  \eta_{ m,n, i} \geq x\Big)  +\ \mathbb{P}\Big(R_{m,n} \geq x, \sum_{i=1}^{n+m}  \eta_{ m,n, i} < x\Big).
\end{align*}
Applying Cram\'{e}r's  moderate deviations for independent random variables (cf. inequality (1) in \cite{FGL13}) to the last inequality,   we deduce that for  all $0\leq x  \leq c \sqrt{ m \wedge n },$
\begin{equation*}
\mathbb{P}\Big(R_{m,n}\geq x\Big)
 \leq  \Big(1- \Phi (x) \Big)\Big(1+ C  \frac{ 1+ x^3   }{  \sqrt{m \wedge n}\  }\Big)   +\ \mathbb{P}\bigg(R_{m,n} \geq x, \sum_{i=1}^{n+m}  \eta_{ m,n, i} < x\bigg) .
\end{equation*}
By Lemma \ref{lama4.9} and  (\ref{norb}), it follows that for all $0 \leq x  \leq  \sqrt{\ln(m \wedge n)} $,
\begin{equation*}
	\mathbb{P}\Big(R_{m,n}\geq x\Big)
 \leq   \Big(1- \Phi (x) \Big)\Big(1+ C  \frac{ 1+ x^3   }{  \sqrt{m \wedge n}\  }\Big).
\end{equation*}
By the last inequality and the inequality $\ln (1+x) \leq  x, x \geq 0,$ we get (\ref{ineq332}) for  all $0 \leq x  \leq  \sqrt{\ln(m \wedge n)} $.

Next, we consider the case  $\sqrt{\ln(m \wedge n)} \leq x  \leq c\, \sqrt{m\wedge n}$.
Clearly,  it holds for all  $x \in \mathbb{R},$
 \begin{align}
\mathbb{P}\bigg(  R_{m,n} \geq x  \bigg)&= \mathbb{P}\bigg(   \sum_{i=1}^{n+m}  \eta_{ m,n, i} + \frac{\ln W_{1, n}}{ n \, V_{m,n,\rho} } - \frac{\ln W_{2,m}}{ m\, V_{m,n,\rho} } \geq x \bigg ) \nonumber \\
 & \leq    I_1+I_2+I_3,   \label{thn635s}
\end{align}
where
 \begin{align*}
 I_1 &= \mathbb{P}\Bigg(  \sum_{i=1}^{n+m}  \eta_{ m,n, i} \geq x\bigg(1- \frac{ (\frac1 n +\frac 1{m \alpha})  x }{\, V_{m,n,\rho}} \bigg)  \Bigg ), \\
 I_2&= \mathbb{P}\Bigg(  \frac{\ln W_{1, n}}{  n \, V_{m,n,\rho} } \geq \frac{   x^{2}}{ n\, V_{m,n,\rho} } \Bigg ) \ \ \ \ \ \ \  \ \textrm{and} \ \ \ \ \ \
  I_3\ =\  \mathbb{P}\Bigg(  -\frac{\ln W_{2,m}}{ m\, V_{m,n,\rho} } \geq \frac{    x^{2}}{m \alpha \, V_{m,n,\rho} } \Bigg ),
\end{align*}
with $\alpha$ given by Lemma \ref{lemma1}.
Next, we give some estimations for  $I_1, I_2$ and $I_3$.  By condition \textbf{A3}, $\sum_{i=1}^{n+m}  \eta_{ m,n, i} $ is a sum of  independent random variables with finite moment generating functions.
 By  Cram\'{e}r's  moderate deviations for independent random variables (cf.\ \cite{FGL13}),   we obtain for all $1\leq x  \leq c\, \sqrt{m\wedge n}$,
 \begin{align*}
I_1 &\leq  \bigg(1- \Phi\Big(x(1- \frac{ (\frac1 n +\frac 1{m \alpha})x }{\, V_{m,n,\rho}} ) \Big)  \bigg )\exp\bigg\{ \frac{  C}{\sqrt{m+n} }\Big(x(1- \frac{ (\frac1 n +\frac 1{m \alpha})x }{\, V_{m,n,\rho}} ) \Big)^{3} \bigg\} \\
&\leq  \bigg(1- \Phi\Big( x(1- \frac{ (\frac1 n +\frac 1{m \alpha})x }{\, V_{m,n,\rho}} )\Big)  \bigg )\exp\bigg\{ C\frac{x^{3}}{\sqrt{m \wedge n} } \bigg\}.
\end{align*}
Using   \eqref{norb},
we deduce  that for all  $x\geq 1$ and $  \varepsilon_n  \in (0, \frac1 2] $,
\begin{align}
  \frac{1-\Phi \left( x (1-  \varepsilon_n) \right)}{1-\Phi \left( x\right) }& = 1+ \frac{ \int_{x (1-  \varepsilon_n)}^x \frac{1}{\sqrt{2\pi}}e^{-t^2/2}dt }{1-\Phi \left( x\right) }\nonumber \\
  &\leq    1+ \frac{\frac{1}{\sqrt{2\pi}} e^{-x^2(1- \varepsilon_n)^2/2}  x \varepsilon_n  }{ \frac{1}{\sqrt{2 \pi} (1+x)} e^{-x^2/2}  }
    \leq   1+  C   x^2    \varepsilon_n   \exp\Big\{ C    x^2   \varepsilon_n  \Big\}   \nonumber \\
    &\leq    \exp\Big\{ 2C    x^2   \varepsilon_n  \Big\}. \label{sfdsh}
\end{align}
Hence, by the fact  $V_{m,n,\rho}\asymp \frac{1}{\sqrt{m}}$, it holds for  all  $1 \leq x  \leq c\, \sqrt{m \wedge n} ,$
\begin{equation}\label{ines4}
I_1  \leq  \Big(1- \Phi(x )  \Big )\exp\Big\{ C \frac{ x^{3} }{\sqrt{m \wedge n} }\Big\} .
\end{equation}
By Markov's inequality, it is easy to see that for all  $x \geq \sqrt{ \ln (m \wedge n)    },$
 \begin{align}
I_2&=\mathbb{P}\Big(  W_{1,n} \geq \exp\big\{    x^{2}  \big \} \Big)\leq  \exp\big\{ -     x^2  \big \}  \mathbb{E} W_{1,n}=\exp\big\{ -     x^2  \big \}  \nonumber  \\
&\leq  C  \frac{ 1+x  }{\sqrt{m \wedge n} } \Big(1-\Phi(x)\Big)   \label{ines5}
\end{align}
and
 \begin{align}
I_3&=\mathbb{P}\Big(  W_{2,m} \leq \exp\big\{ -  \alpha^{-1}  x^{2}  \big \} \Big) \leq  \exp\big\{ -     x^2  \big \}  \mathbb{E} W_{2,m}^{-\alpha}\leq  C \exp\big\{ -    x^2  \big \}  \nonumber  \\
&\leq  C  \frac{ 1+x  }{\sqrt{m \wedge n} } \Big(1-\Phi(x)\Big).   \label{ines6}
\end{align}
Combining (\ref{ines4})-(\ref{ines6}) together, we obtain for  all  $\sqrt{ \ln (m \wedge n) }\leq  x \leq c\, \sqrt{m \wedge n} ,$
 \begin{align*}
\mathbb{P}\bigg( R_{m,n}   \geq x  \bigg )
&\leq \Big(1- \Phi(x )  \Big )\exp\Big\{ C_1 \frac{ x^{3} }{\sqrt{m \wedge n} }\Big\} + C_2  \frac{(1+x) }{\sqrt{m \wedge n} } \Big(1-\Phi(x)\Big)\\
 &\leq \Big(1- \Phi(x )  \Big )\exp\Big\{ C_3 \frac{ x^{3} }{\sqrt{m \wedge n} }\Big\},
\end{align*}
which implies the desired inequality for  all $\sqrt{ \ln (m \wedge n) }\leq  x \leq c\, \sqrt{m \wedge n} .$    \hfill\qed

The following lemma gives the lower bound in Theorem \ref{th2.2}.
\begin{lemma}\label{th100}
Assume that the conditions \textbf{A3} and \textbf{A4} are satisfied. Then  it holds
 for   all  $0 \leq  x \leq  c \, \sqrt{m \wedge n}  ,$
\begin{equation}\label{ineq3d}
\ln \frac{\mathbb{P}\big(   R_{m,n} \geq x  \big)}{1-\Phi(x)}  \geq - C  \frac{1+x^3  }{\sqrt{ m \wedge  n }}.
\end{equation}
\end{lemma}
\emph{Proof.}  The proof for the lower bound is similar  to that  of upper bound.
For instance, to prove (\ref{ineq3d}) for all  $\sqrt{ \ln (m \wedge n) } \leq  x \leq  c \, \sqrt{m \wedge n}  ,$   we only need to notice that
 \begin{align*}
\mathbb{P}\bigg(  R_{m,n} \geq x  \bigg)
&=  \mathbb{P}\bigg(   \sum_{i=1}^{n+m}  \eta_{ m,n, i} + \frac{\ln W_{1, n}}{ n\,V_{m,n,\rho} } - \frac{\ln W_{2,m}}{ m\, V_{m,n,\rho} } \geq x \bigg )  \\
  &\geq    I_4-I_5-I_6,
\end{align*}
where
 \begin{align*}
 I_4 &= \mathbb{P}\Bigg(  \sum_{i=1}^{n+m}  \eta_{ m,n, i} \geq x\bigg(1+ \frac{ (\frac1 {n\alpha} +\frac 1{m})  x }{\,V_{m,n,\rho}} \bigg)  \Bigg ), \\
 I_5&= \mathbb{P}\Bigg( - \frac{\ln W_{1, n}}{  n  \,V_{m,n,\rho} } \geq  \frac{   x^{2}}{ n\alpha \,V_{m,n,\rho} } \Bigg ) \quad \textrm{and} \quad
  I_6 =  \mathbb{P}\Bigg(   \frac{\ln W_{2,m}}{ m  \,V_{m,n,\rho} } \geq  \frac{    x^{2}}{m \,V_{m,n,\rho} } \Bigg ),
\end{align*}
with $\alpha$ given by Lemma \ref{lemma1}. For the  remaining of the proof,
we can use  the argument similar to  the proof of Lemma  \ref{th200}.
 \hfill\qed

\section*{Conflict of interest}
The authors declared that they have no conflict of interest.

\section*{Acknowledgements}
Fan would like to thank professor Quansheng Liu for his helpful discussion on the harmonic moments for branching processes in a random environment.
This work has been partially supported by the National Natural Science Foundation
of China (Grant no.\ 11971063).

\end{document}